\documentclass{amsart}

\usepackage{a4wide}
\usepackage{amsmath,amscd}%
\usepackage{amsfonts}%
\usepackage{amssymb}%
\usepackage{graphicx}
\usepackage[usenames,dvipsnames]{color}
\usepackage{subfigure}
\usepackage{enumerate}
\usepackage[noprefix]{nomencl}  
\makenomenclature 
\usepackage[heads=vee]{diagrams}
\usepackage{hyperref}
\hypersetup{colorlinks,citecolor=black,filecolor=black,linkcolor=black,urlcolor=black}

\newcommand{\sfa}{\mathsf{a}}
\newcommand{\sfb}{\mathsf{b}}
\newcommand{\sfc}{\mathsf{c}}
\newcommand{\sfe}{\mathsf{e}}
\newcommand{\sff}{\mathsf{f}}
\newcommand{\sfg}{\mathsf{g}}
\newcommand{\sfh}{\mathsf{h}}

\newcommand{\kdifform}[2]{#1^{(#2)}} 
\newcommand{\difform}[1]{#1} 
\newcommand{\kchain}[2]{\mathbf{#1}_{(#2)}} 
\newcommand{\kcochain}[2]{\mathbf{#1}^{(#2)}} 
\newcommand{\reconstruction}{\mathcal{I}} 
\newcommand{\reduction}{\mathcal{R}} 
\newcommand{\kformspace}[1]{\Lambda^{#1}} 
\newcommand{\kformspacedomain}[2]{\Lambda^{#1}(#2)} 
\newcommand{\kchainspacedomain}[2]{C_{#1}(#2)} 
\newcommand{\kcochainspacedomain}[2]{C^{#1}(#2)} 
\newcommand{\projection}{\pi_{h}} 
\newcommand{\pullback}{\Phi^{\star}} 

\newcommand{\ederiv}{\mathrm{d}} 
\newcommand{\ud}{\mathrm{d}} 
\newcommand{\p}{\partial}
\newcommand{\spacemap}{\rightarrow} 
 
\newcommand{\duality}[2]{\langle #1, #2\rangle} 

\renewcommand{\eqref}[1]{(\ref{#1})} 
\newcommand{\figref}[1]{Figure~\ref{#1}} 
\newcommand{\secref}[1]{Section~\ref{#1}} 
\newcommand{\defref}[1]{Definition~\ref{#1}} 
\newcommand{\propref}[1]{Proposition~\ref{#1}} 
\newcommand{\lemmaref}[1]{Lemma~\ref{#1}} 

\newtheorem{theorem}{Theorem}
\theoremstyle{plain}

\newtheorem{corollary}{Corollary}

\newtheorem{definition}{Definition}

\newtheorem{lemma}{Lemma}

\newtheorem{proposition}{Proposition}
\newtheorem{remark}{Remark}

\numberwithin{equation}{section}

\newcommand{\define}{\mathrel{\mathop:}=}

\newcommand{\norm}[1]{\Vert #1\Vert}
\newcommand{\tr}{\mathrm{tr\;}}

\newcommand{\ve}{\varepsilon}

\graphicspath{{figures/}}

\begin{document}
\title[A priori error estimates for compatible spectral discretization of the Stokes problem]{A priori error estimates for compatible spectral discretization of the Stokes problem for all admissible boundary conditions}
\author{Jasper Kreeft}
\address[]
	{Delft University of Technology, Faculty of Aerospace Engineering, \newline%
	\indent  Kluyverweg 2, 2629 HT Delft, The Netherlands.}%
\email[]{J.J.kreeft@gmail.com, M.I.Gerritsma@TUDelft.nl}%
\author{Marc Gerritsma}
\thanks{Jasper Kreeft is funded by STW Grant 10113}
\thanks{This paper is in final form and no version of it will be submitted for
publication elsewhere.}
\date{\today}
\subjclass{Primary 76D07, 65N30; Secondary 65M70, 12Y05, 13P20} %
\keywords{Stokes problem, mixed finite elements, mimetic/compatible discretization, error estimates}%

\begin{abstract}
This paper describes the recently developed mixed mimetic spectral element method for the Stokes problem in the vorticity-velocity-pressure formulation. This compatible discretization method relies on the construction of a conforming discrete Hodge decomposition, that is based on a bounded projection operator that commutes with the exterior derivative. The projection operator is the composition of a reduction and a reconstruction step. The reconstruction in terms of mimetic spectral element basis-functions are tensor-based constructions and therefore hold for curvilinear quadrilateral and hexahedral meshes.

For compatible discretization methods that contain a conforming discrete Hodge decomposition, we derive optimal a priori error estimates which are valid for all admissible boundary conditions on both Cartesian and curvilinear meshes. These theoretical results are confirmed by numerical experiments. These clearly show that the mimetic spectral elements outperform the commonly used $H(\mathrm{div})$-compatible Raviart-Thomas elements.
\end{abstract}

\maketitle

\section{Introduction}
Let $\Omega\subset\mathbb{R}^n$, $n\geq2$, be a bounded contractible domain with boundary $\Gamma=\p\Omega$. On this domain we consider the Stokes problem, consisting of the equations for conservation of momentum and for conservation of mass,
\begin{subequations}
\begin{align}
\nabla\cdot\sigma&=\vec{f}\quad\mathrm{on}\ \Omega,\\
\mathrm{div}\,\vec{u}&=g\quad\mathrm{on}\ \Omega,
\end{align}
\end{subequations}
where the stress tensor $\sigma$ is given by
\begin{equation}
\sigma=-\nu\nabla\vec{u}+pI,
\end{equation}
with $\vec{u}$ the velocity vector, $p$ the pressure, $\vec{f}$ the forcing term, $g$ the mass source and $\nu$ the kinematic viscosity. For analysis purposes we choose $\nu=1$.

This paper considers the recently developed mixed mimetic spectral element method (MMSEM) \cite{kreeft2012,kreeftpalhagerritsma2011}. This compatible finite/spectral element method is based on the compatible discretization of the exterior derivative $\ud$ from differential geometry, which represents the vector operators, grad, curl and div. The Stokes problem expressed in terms of these vector operations is known as the \emph{vorticity-velocity-pressure} (VVP) formulation, \cite{bochevgunzburger2009,dubois2002}. For the VVP formulation, the Laplace operator is split using the vector identity, $-\Delta\vec{u}=\mathrm{curl}\,\mathrm{curl}^*\,\vec{u}-\mathrm{grad}^*\,\mathrm{div}\,\vec{u}$, and by introducing vorticity as auxiliary variable, $\vec{\omega}=\mathrm{curl}^*\,\vec{u}$. The VVP formulation of the Stokes problem becomes
\begin{subequations}
\label{stokessinglevector}
\begin{align}
\vec{\omega}-\mathrm{curl}^*\,\vec{u}&=0,\quad\mathrm{on}\ \Omega\label{stokessinglevector1}\\
\mathrm{curl}\,\vec{\omega}-\mathrm{grad}^*\,\mathrm{div}\,\vec{u}+\mathrm{grad}^*\,p&=\vec{f},\quad\mathrm{on}\ \Omega \label{stokessinglevector2}\\ 
\mathrm{div}\,\vec{u}&=g,\quad\mathrm{on}\ \Omega.\label{stokessinglevector3}
\end{align}
\end{subequations}
Following \cite{bochevgunzburger2009,kreeft2012} we make a distinction between the operators grad, curl and div, that correspond to the classical Newton-Leibnitz, Stokes circulation and Gauss divergence theorems, and the operators -grad$^*$, curl$^*$ and -div$^*$ that are their formal Hilbert adjoints,
\[
\big(\vec{a},-\mathrm{grad}^*\,b\big):=\big(\mathrm{div}\,\vec{a},b\big),\quad \big(\vec{a},\mathrm{curl}^*\,\vec{b}\big):=\big(\mathrm{curl}\,\vec{a},\vec{b}\big),\quad \big(a,-\mathrm{div}^*\,\vec{b}\big):=\big(\mathrm{grad}\,a,\vec{b}\big).
\]
The distinction between the two types of differential operators is made explicitly, because the construction of our conforming finite element spaces relies on the three mentioned integration theorems, while the mixed formulation relies on the formal Hilbert adjoint relations. While in vector calculus this distinction is not common, in differential geometry these structures naturally appear since they make a clear distinction between metric-free (topological) and metric-dependent operations.

The MMSEM is a compatible discretization method that relies on the construction of a conforming discrete Hodge-decomposition, which implies a discrete Poincar\'e inequality. It requires the development of a bounded projection operator that commutes with the exterior derivative. The bounded projection is a composition of a reduction by means of integration and mimetic spectral element basisfunctions as reconstruction.

The reduction onto $k$-dimensional submanifolds result in the discrete unknowns representing \emph{integral quantities}. This is one of the major differences with related methods as the Marker and Cell scheme \cite{harlowwelch1965} and the lowest-order Raviart-Thomas and N\'ed\'elec compatible finite elements \cite{nedelec1980,raviartthomas1977}, where use is made of averaged quantities.

The basis functions, used for the reconstruction, are constructed using tensor products of one dimensional nodal and edge interpolation basis functions \cite{gerritsma2011}, and therefore hold for quadrilateral and hexahedral meshes. They belong to the class of compatible finite elements, and were constructed based on the mimetic framework first described in \cite{hymanscovel1988} and later extended in \cite{bochevhyman2006}. The mimetic framework, including the mimetic spectral elements, were extensively described in \cite{kreeftpalhagerritsma2011}. This mimetic framework relies on the languages of differential geometry instead of vector calculus, and algebraic topology as its discrete counterpart. 

The use of differential geometry and algebraic topology enjoys increasing popularity for the development of compatible schemes, \cite{arnoldfalkwinther2006,arnoldfalkwinther2010,bochevhyman2006,bossavit1998,bossavit9900,desbrun2005c,hiptmair2001,hiptmair2002}. Compatible discretizations are often combined with mixed formulations. Mixed formulations are described extensively in among others \cite{brezzifortin,giraultraviart} and in terms of differential forms in \cite{arnoldfalkwinther2006,arnoldfalkwinther2010} for the Hodge-Laplacian and in \cite{kreeft2012} for the VVP formulation of the Stokes problem.


The MMSEM contains compatible finite elements that are compatible with all admissible types of boundary conditions for the Stokes problem in VVP formulation. We will show that the method obtains optimal rates of convergence for all variables on curvilinear meshes and for all admissible boundary conditions, i.e. standard and nonstandard. It is therefore extending the error estimates found in literature, which are often specifically constructed for certain types of boundary conditions, \cite{abboud2011,arnold2011,bochevgunzburger,bramble1994,dubois2003b,girault1988}. To show optimal convergence a priori error estimates are derived.

This is an improvement with respect to the well-known Raviart-Thomas compatible finite elements. These are not compatible in case of Dirichlet boundary conditions and therefore lead to suboptimal convergence behavior, as was shown in \cite{arnold2011,dubois2003b}. This non-compatibility results in a decrease in rate of convergence of maximal $\tfrac{3}{2}$ order.

From a physical/fluid dynamics point-of-view the new method is relevant because it combines optimal convergence with a pointwise divergence-free discretization (in absence of any mass source) of arbitrary order on curvilinear meshes, valid for all allowable types of boundary conditions, among which the no-slip condition.

The derived rates of convergence are confirmed using simple manufactured solution problems, discretized on both Cartesian and curvilinear meshes. The fact that the analysis holds for all admissible boundary conditions is also reflected in the numerical results.

This paper is organized as follows: First an introduction into differential geometry is given and the Stokes problem is reformulated in terms of differential forms. In \secref{sec:mixedformulation} the mixed formulation is given and well-posedness is proven. In \secref{sec:discretization} the key properties of the mimetic discretization are explained that lead to compatible function spaces. This includes a discussion on the relevant properties of algebraic topology, the definitions op mimetic operators, the introduction of mimetic spectral element basisfunctions and finally the proof of discrete well-posedness. Having formulated the conforming/compatible finite element spaces, the error estimates are developed in \secref{sec:errorestimates} and the numerical results are shown in \secref{sec:numericalresults}.

\section{Notation and preliminaries}\label{sec:differentialgeometry}

\subsection{Differential forms}
Differential forms offer significant benefits in the construction of structure-preserving spatial discretizations. For example, the coordinate-free action of the exterior derivative and generalized Stokes theorem give rise to commuting properties with respect to mappings between different manifolds. Acknowledging and respecting these kind of commuting properties are essential for the structure preserving behavior of the mimetic method.

Only those concepts from differential geometry which play a role in the remainder of this paper will be explained. More can be found in \cite{abrahammarsdenratiu,flanders,frankel,kreeftpalhagerritsma2011}.

Let $\Lambda^k(\Omega)$ denote a space of \emph{differential $k$-forms} or \emph{$k$-forms}, on a sufficiently smooth bounded $n$-dimensional oriented manifold $\Omega\subset\mathbb{R}^n$ with boundary $\Gamma=\p\Omega$. Every element $a\in\Lambda^k(\Omega)$ has a unique representation of the form
\begin{equation}
\label{differentialform}
a=\sum_If_I(\mathbf{x})\ud x^{i_1}\wedge\ud x^{i_2}\wedge\cdots\wedge\ud x^{i_k},
\end{equation}
where $I=i_1,\hdots,i_k$ with $1\leq i_1<\hdots<i_k\leq n$ and where $f_I(\mathbf{x})$ is a continuously differentiable scalar function, $f_I(\mathbf{x})\in \mathcal{C}^\infty(\Omega)$. Differential $k$-forms are naturally integrated over $k$-dimensional manifolds, i.e. for $a\in \Lambda^k(\Omega)$ and $\Omega_k\subset\mathbb{R}^n$, with $k=\mathrm{dim}(\Omega_k)$,
\begin{equation}
\label{integration}
\int_{\Omega_k}a\in\mathbb{R}\quad\Leftrightarrow\quad\langle a,\Omega_k\rangle\in\mathbb{R},
\end{equation}
where $\langle\cdot,\cdot\rangle$ indicates a duality pairing between the differential form and the geometry. Note that the $n$-dimensional computational domain is indicated as $\Omega$, so 	without subscript. The differential forms live on manifolds and transform under the action of mappings. Let $\Phi:\widehat{\Omega}\rightarrow\Omega$ be a mapping between two manifolds. Then we can define the pullback operator, $\Phi^\star:\Lambda^k(\Omega)\rightarrow\Lambda^k(\widehat{\Omega})$, expressing the $k$-form on the $n$-dimensional reference manifold, $\widehat{\Omega}$. The mapping, $\Phi$, and the pullback, $\Phi^\star$, are each others formal adjoints with respect to a duality pairing \eqref{integration},
\begin{equation}
\int_{\Phi(\widehat{\Omega}_l)}a=\int_{\widehat{\Omega}_l}\Phi^\star a\quad\Leftrightarrow\quad\langle a,\Phi(\widehat{\Omega}_l)\rangle=\langle\Phi^\star a,\widehat{\Omega}_l\rangle,
\end{equation}
where $\widehat{\Omega}_l$ is an $l$-dimensional submanifold of $\widehat{\Omega}$ and $\Omega_k=\Phi(\widehat{\Omega}_l)$ a $k$-dimensional submanifold of $\Omega$. A special case of the pullback operator is the trace operator. The trace of $k$-forms to the boundary, $\tr:\Lambda^k(\Omega)\rightarrow\Lambda^k(\p\Omega)$, is the pullback of the inclusion of the boundary of a manifold, $\p\Omega\hookrightarrow\Omega$, see \cite{kreeftpalhagerritsma2011}.

The wedge product, $\wedge$, of two differential forms $a\in\Lambda^k(\Omega)$ and $b\in\Lambda^l(\Omega)$ is a mapping: $\wedge:\Lambda^k(\Omega)\times\Lambda^l(\Omega)\rightarrow\Lambda^{k+l}(\Omega),\ k+l\leq n$. The wedge product is a skew-symmetric operator, i.e. $a\wedge b=(-1)^{kl}b\wedge a$.

An important operator in differential geometry is the exterior derivative, $\ud:\Lambda^k(\Omega)\rightarrow\Lambda^{k+1}(\Omega)$. It is induced by the \emph{generalized Stokes' theorem}, combining the classical Newton-Leibnitz, Stokes circulation and Gauss divergence theorems. Let $\Omega_{k+1}$ be a $(k+1)$-dimensional manifold and $a\in\Lambda^k(\Omega)$, then
\begin{equation}\label{stokestheorem}
\int_{\p\Omega_{k+1}}a=\int_{\Omega_{k+1}}\ud a\quad\Leftrightarrow\quad \langle a,\partial\Omega_{k+1}\rangle=\langle\ud a,\Omega_{k+1}\rangle,
\end{equation}
where $\partial\Omega_{k+1}$ is a $k$-dimensional manifold being the boundary of $\Omega_{k+1}$. The duality pairing in \eqref{stokestheorem} shows that the exterior derivative is the formal adjoint of the \emph{boundary operator} $\p:\Omega_{k+1}\rightarrow\Omega_k$. The exterior derivative is independent of any metric and coordinate system. Applying the exterior derivative twice always leads to the null $(k+2)$-form, $\ud(\ud a)=0^{(k+2)}$ for all $a\in\Lambda^k(\Omega)$.  As a consequence, on contractible domains the exterior derivative gives rise to an exact sequence, called \emph{De Rham complex} \cite{frankel}, and indicated by $(\Lambda,\ud)$,
\begin{equation}
\mathbb{R}\hookrightarrow\Lambda^0(\Omega)\stackrel{\ud}{\longrightarrow}\Lambda^1(\Omega)\stackrel{\ud}{\longrightarrow}\cdots\stackrel{\ud}{\longrightarrow}\Lambda^n(\Omega)\stackrel{\ud}{\longrightarrow}0.
\label{derhamcomplex}
\end{equation}
In vector calculus a similar sequence exists, where, from left to right for $\mathbb{R}^3$, the $\ud$'s denote the vector operators grad, curl and div. The exterior derivative and wedge product are related according to Leibnitz's rule as: for all $a\in \Lambda^k(\Omega)$ and $b\in\Lambda^l(\Omega)$,
\begin{equation}
\ud\big(a\wedge b\big)=\ud a\wedge b+(-1)^ka\wedge\ud b,\quad\mathrm{for}\ k+l<n.
\label{leibnitz}
\end{equation}
The pullback operator and exterior derivative possess the following commuting property,
\begin{equation}
\Phi^\star\ud a=\ud\Phi^\star a,\quad\forall a\in\Lambda^k(\Omega).
\end{equation}

In this paper we will consider Hilbert spaces $L^2\Lambda^k(\Omega)\supset\Lambda^k(\Omega)$, where in \eqref{differentialform} the functions $f_I(\mathbf{x})\in L^2(\Omega)$. The pointwise inner-product of $k$-forms, $(\cdot,\cdot)$, is constructed using inner products of one-forms, that is based on the inner product on vector spaces, see \cite{flanders,frankel}. The wedge product and inner product induce the Hodge-$\star$ operator, $\star:L^2\Lambda^k(\Omega)\rightarrow L^2\Lambda^{n-k}(\Omega)$, a metric operator that includes orientation. Let $a,b\in L^2\Lambda^k(\Omega)$, then
\begin{equation}
a\wedge\star b:=\big(a,b\big)\sigma,
\label{hodgestar}
\end{equation}
where $\sigma\in\Lambda^n(\Omega)$ is a unit volume form, $\sigma=\star1$. In geometric physics the Hodge-$\star$ switches between an inner-oriented description of physical variables and an outer-oriented description. See \cite{kreeft2012,kreeftpalhagerritsma2011,mattiussi2000,tonti1} for a thorough discussion on the concepts of inner and outer orientation. The space of square integrable $k$-forms on $\Omega$ can be equipped with a $L^2$ inner product, $\big(\cdot,\cdot\big)_\Omega:L^2\Lambda^k(\Omega)\times L^2\Lambda^k(\Omega)\rightarrow\mathbb{R}$, given by,
\begin{equation}
\big(a,b\big)_\Omega:=\int_\Omega\big(a,b\big)\kdifform{\sigma}{n}=\int_\Omega a\wedge\star b.
\label{L2innerproduct}
\end{equation}
The norm corresponding to the space $L^2\Lambda^k(\Omega)$ is $\norm{a}_{L^2\Lambda^k}=\sqrt{\big(a,a\big)_\Omega}$. Higher degree Sobolev spaces, $H^m\Lambda^k$, consists of all $k$-forms as in \eqref{differentialform} where $f_I(\mathbf{x})\in H^m(\Omega)$, with corresponding norms $|a|_{H^m\Lambda^k}$ and $\Vert a\Vert_{H^m\Lambda^k}$. The Hilbert space associated to the exterior derivative $H\Lambda^k(\Omega)$ is defined as
\begin{equation}
H\Lambda^k(\Omega)=\{a\in L^2\Lambda^k(\Omega)\;|\;\ederiv a\in L^2\Lambda^{k+1}(\Omega)\}.
\end{equation}
and the norm corresponding to $H\Lambda^k(\Omega)$ is defined as $\Vert a\Vert^2_{H\Lambda^k}:=\Vert a\Vert^2_{L^2\Lambda^k}+\Vert\ederiv a\Vert^2_{L^2\Lambda^{k+1}}$. The $H\Lambda^k$-semi-norm is the $L^2$-norm of the exterior derivative, $\vert a\vert_{H\Lambda^k}=\Vert\ederiv a\Vert_{L^2\Lambda^{k+1}}$. Note that $H^1\Lambda^k(\Omega)\subseteq H\Lambda^k(\Omega)\subseteq L^2\Lambda^k(\Omega)$, where the left equality holds for $k=0$ and the right for $k=n$. The $L^2$-de Rham complex, also called \emph{Hilbert complex} \cite{bruninglesch1992}, $(H\Lambda,\ud)$, is the exact sequence of maps and spaces given by
\begin{equation}
\mathbb{R}\hookrightarrow H\Lambda^0(\Omega)\stackrel{\ederiv}{\longrightarrow} H\Lambda^1(\Omega)\stackrel{\ederiv}{\longrightarrow}\cdots\stackrel{\ederiv}{\longrightarrow} H\Lambda^n(\Omega)\stackrel{\ud}{\longrightarrow}0.
\end{equation}
In terms of vector operations the Hilbert complex becomes for $\Omega\subset\mathbb{R}^3$,
\[
H^1(\Omega)\stackrel{\rm grad}{\longrightarrow} H(\mathrm{curl},\Omega)\stackrel{\rm curl}{\longrightarrow}H(\mathrm{div},\Omega)\stackrel{\rm div}{\longrightarrow} L^2(\Omega),
\]
and for $\Omega\subset\mathbb{R}^2$, either
\[
H^1(\Omega)\stackrel{\rm grad}{\longrightarrow} H(\mathrm{rot},\Omega)\stackrel{\rm rot}{\longrightarrow}L^2(\Omega),\quad\mathrm{or}\quad
H^1(\Omega)\stackrel{\rm curl}{\longrightarrow} H(\mathrm{curl},\Omega)\stackrel{\rm div}{\longrightarrow}L^2(\Omega).
\]
The two are related by the Hodge-$\star$ operator \eqref{hodgestar}, see \cite{palha2010},
\begin{equation}
\label{doublehilbertcomplex}
\!\!\!\!\!\!\!
\begin{matrix}
H\Lambda^{0}(\Omega)\!\!&\!\!\stackrel{\ederiv}{\longrightarrow}\!\!&\!\! H\Lambda^{1}(\Omega)\!\!&\!\!\stackrel{\ederiv}{\longrightarrow}\!\!&\!\!L^2\Lambda^{2}(\Omega)\\
\star\updownarrow & & \star\updownarrow & &\star\updownarrow   \\
L^2\Lambda^{2}(\Omega)\!\!&\!\!\stackrel{\ederiv}{\longleftarrow}\!\!&\!\!H\Lambda^{1}(\Omega)\!\!&\!\!\stackrel{\ederiv}{\longleftarrow}\!\!&\!\!H\Lambda^{0}(\Omega)
\end{matrix}
\quad\Leftrightarrow\quad
\begin{matrix}
H^1(\Omega)\!\!&\!\!\stackrel{\mathrm{curl}}{\longrightarrow}\!\!&\!\!H(\mathrm{curl},\Omega)\!\!&\!\!\stackrel{\mathrm{div}}{\longrightarrow}\!\!&\!\!L^2(\Omega)\\
\star\updownarrow & & \star\updownarrow & &\star\updownarrow   \\
L^2(\Omega)\!\!&\!\!\stackrel{\mathrm{rot}}{\longleftarrow}\!\!&\!\!H(\mathrm{rot},\Omega)\!\!&\!\!\stackrel{\mathrm{grad}}{\longleftarrow}\!\!&\!\!H^1(\Omega).
\end{matrix}
\end{equation}
\begin{remark}
The upper complex is associated with outer-oriented $k$-forms, i.e. $k$-forms that are associated with outer-oriented manifolds, and the lower complex is associated with inner-oriented $k$-forms. 
In this paper we mainly consider the upper complex and circumvent the lower complex by means of integration by parts. Only the pressure and tangential velocity boundary conditions are given on the lower complex, as we will see in the following sections.
\end{remark}
A similar double Hilbert complex can be constructed in $\mathbb{R}^3$. Since the exterior derivative is nilpotent, it ensures that the range, $\mathcal{B}^k:=\ud\, H\Lambda^{k-1}(\Omega)$, of the exterior derivative on $(k-1)$-forms is contained in the nullspace, $\mathcal{Z}^k:=\{\;a\in H\Lambda^k(\Omega)\;|\;\ud a=0\;\}$, of the exterior derivative on $k$-forms, $\mathcal{B}^k\subseteq\mathcal{Z}^k$.

Every space of $k$-forms in the complex $(H\Lambda,\ud)$ can be decomposed into the nullspace of $\ud$, and its orthogonal complement, $H\Lambda^k(\Omega)=\mathcal{Z}^k\oplus\mathcal{Z}^{k,\perp}$. This is the \emph{Hodge decomposition}, where on contractible domains $\mathcal{Z}^k=\mathcal{B}^k$. By the Hodge decomposition it follows that the exterior derivative is an isomorphism $\ud:\mathcal{Z}^{k,\perp}\rightarrow\mathcal{B}^{k+1}$.

The inner product gives rise to the formal Hilbert adjoint of the exterior derivative, the codifferential operator, $\ud^*:H^*\Lambda^k(\Omega)\rightarrow L^2\Lambda^{k-1}(\Omega)$. Let $H^*\Lambda^k(\Omega)=\{\,a\in L^2\Lambda^k(\Omega)\,|\,\ud^*a\in L^2\Lambda^{k-1}(\Omega)\}$, then
\begin{equation}
\big(\ud a,b\big)_\Omega=(a,\ud^* b)_\Omega,\quad\forall\, a\in H\Lambda^{k-1}(\Omega),\ b\in H^*\Lambda^k.
\end{equation}
In case of non-zero trace, and by combining \eqref{L2innerproduct}, \eqref{stokestheorem} and \eqref{leibnitz}, we obtain integration by parts,
\begin{equation}
\big(a,\ud^*b\big)_\Omega=\big(\ud a,b\big)_\Omega-\int_{\p\Omega} \tr a\wedge \tr\star b.
\label{integrationbyparts}
\end{equation}
Also the codifferential operator is nilpotent, $\ud^*(\ud^*\kdifform{a}{k})=0$, i.e., its range is contained in its nullspace, $\mathcal{B}^{*,k}\subseteq\mathcal{Z}^{*,k}$, where $\mathcal{B}^{*,k}:=\ud^*H^*\Lambda^{k+1}(\Omega)$ and $\mathcal{Z}^{*,k}:=\{\; a\in H^*\Lambda^k(\Omega)\;|\;\ud^*a=0\;\}$. In fact the codifferential is an isomorphism $\ud^*:\mathcal{Z}^{*,k,\perp}\rightarrow\mathcal{B}^{*,k+1}$, where $\mathcal{Z}^{*,k,\perp}$ follows from the following Hodge decomposition, $\Lambda^k(\Omega)=\mathcal{Z}^{*,k}\oplus\mathcal{Z}^{*,k,\perp}$. On contractible manifolds this gives rise to the following exact sequence,
\begin{equation}
0\stackrel{\ud^*}{\longleftarrow}H^*\Lambda^0(\Omega)\stackrel{\ud^*}{\longleftarrow}H^*\Lambda^1(\Omega)\stackrel{\ud^*}{\longleftarrow}\cdots\stackrel{\ud^*}{\longleftarrow}H^*\Lambda^n(\Omega)\hookleftarrow\mathbb{R}.
\end{equation}
In vector notation from right to left the $\ud^*$'s denote the -grad$^*$, curl$^*$ and -div$^*$ operators in $\mathbb{R}^3$, as were also mentioned in the introduction. However, whereas the exterior derivative is a metric-free operator, the codifferential operator is metric-dependent. The Hodge-Laplace operator, $\Delta:H^2\Lambda^k(\Omega)\rightarrow L^2\Lambda^k(\Omega)$, is constructed as a composition of the exterior derivative and the codifferential operator,
\begin{equation}
\label{laplace}
-\Delta\, a:=(\ud^*\ud+\ud\ud^*)\,a.
\end{equation}
An important inequality in stability analysis, relating the $L^2\Lambda^k$-norm and the $H\Lambda^k$-norm, is Poincar\'e inequality.
\begin{lemma}[\textbf{Poincar\'e inequality}]\label{lem:poicareinequality}\cite{arnoldfalkwinther2010}
Consider the Hilbert complex $(H\Lambda,\ederiv)$, then the exterior derivative is a bounded bijection from $\mathcal{Z}^{k,\perp}$ to $\mathcal{B}^{k+1}$, and hence, by Banach's bounded inverse theorem, there exists a constant $c_P$ such that
\begin{equation}
\Vert a\Vert_{H\Lambda^k}\leq c_P\Vert\ederiv a\Vert_{L^2\Lambda^{k+1}},\quad\forall a\in \mathcal{Z}^{k,\perp}.
\end{equation}
\end{lemma}
Finally, for Hilbert spaces with essential boundary conditions we write, $H_0\Lambda^k(\Omega):=\{\,a\in H\Lambda^k(\Omega)\,|\,\tr a=0\,\}$, and for natural boundary conditions we consider the following trace map, $\tr\star:H\Lambda^k(\Omega)\rightarrow H^{\frac{1}{2}}\Lambda^{n-k}(\p\Omega)$.

\subsection{Stokes problem in differential form notation}
Consider again a bounded contractible domain $\Omega\subset\mathbb{R}^n$. Because we require exact conservation of mass and because we can perform exact discretization of the exterior derivative, see \secref{mimeticoperators}, we use the following formulation for the Stokes problem: let $(\omega,u,p)\in\{\Lambda^{n-2}(\Omega)\times\Lambda^{n-1}(\Omega)\times\Lambda^{n}(\Omega)\}$, then the VVP formulation is given by
\begin{subequations}
\label{stokeseq}
\begin{align}
\omega-\ud^* u&=0,\quad\mathrm{in}\ \Lambda^{n-2}(\Omega),\label{stokeseq1}\\
\ud^*\ud u+\ud\omega+\ud^*p&=f,\quad\mathrm{in}\ \Lambda^{n-1}(\Omega),\label{stokeseq2}\\
\ud u&=g,\quad\mathrm{in}\ \Lambda^n(\Omega).\label{stokeseq3}
\end{align}
\end{subequations}
In the VVP formulation the pressure in \eqref{stokeseq2} acts as a Lagrange multiplier for the constraint on velocity, \eqref{stokeseq3}, whereas velocity in \eqref{stokeseq1} acts as a Lagrange multiplier for the constraint on vorticity in \eqref{stokeseq2}.

Let $\Gamma=\p\Omega$ be the boundary of $\Omega$, where
\[
\Gamma=\Gamma_\omega\cup\Gamma_t,\quad \Gamma_\omega\cap\Gamma_t=\emptyset,\quad\mathrm{and}\quad
\Gamma=\Gamma_n\cup\Gamma_\pi,\quad\Gamma_n\cap\Gamma_\pi=\emptyset.
\]
We will impose the tangential vorticity and normal velocity as essential boundary conditions, and the tangential velocity and the pressure plus divergence of velocity as the natural boundary conditions:
\begin{subequations}
\label{boundaryconditions}
\begin{align}
\tr\omega&=0\quad\quad\  \mathrm{on}\ \Gamma_\omega,\\
\tr u&=0\quad\quad\  \mathrm{on}\ \Gamma_n,\\
\quad\tr\star u&=u_{b,t}\quad\ \mathrm{on}\ \Gamma_t,\quad\; \mathrm{with}\ u_{b,t}\in H^{\frac{1}{2}}\Lambda^{1}(\Gamma_t),\\
\tr\star(\ud u+p)&=\Pi_b\quad\ \ \, \mathrm{on}\ \Gamma_\pi,\quad \mathrm{with}\ \Pi_b\in H^{\frac{1}{2}}\Lambda^0(\Gamma_\pi).
\end{align}
\end{subequations}
Then the boundary $\Gamma$ can be partitioned into four sections, $\Gamma=\bigcup_{i=1}^4\Gamma_i$, with $\Gamma_i\cap\Gamma_j=\emptyset$ for $i\neq j$, where
\begin{equation}
\label{boundarydecomposition}
\Gamma_1:=\Gamma_t=\Gamma_n,\quad\Gamma_2:=\Gamma_t=\Gamma_\pi,\quad\Gamma_3:=\Gamma_\omega=\Gamma_n,\quad\Gamma_4:=\Gamma_\omega=\Gamma_\pi.
\end{equation}
This decomposition, introduced before in \cite{dubois2002,hughesfranca1987,kreeft2012}, shows all admissible boundary conditions. It will also follow directly from the mixed formulation, see \eqref{mixedstokes}, Section~\ref{sec:mixedformulation}.

In case, $\Gamma=\Gamma_1\cap\Gamma_3$, $\Gamma_2\cup\Gamma_4=\emptyset$, no pressure boundary conditions are prescribed, and so the pressure is only determined up to an element $\hat{p}\in\mathcal{Z}^{*,n}$, i.e. up to a constant. As a post processing step either the pressure in a point in $\Omega$ can be set, or a zero average pressure can be imposed; i.e. $\int_\Omega \hat{p}=0$. In case $\Gamma=\Gamma_4$, no velocity boundary conditions are prescribed, and so the solution of velocity is determined modulo a curl$^*$-free element, i.e. modulo $\hat{u}\in\mathcal{Z}^{*,n-1}$.

\section{Mixed formulation}\label{sec:mixedformulation}

\subsection{Mixed formulation of Stokes problem}
The use of a mixed formulation is based on the following reasoning; We know how to discretize exactly the metric-free exterior derivative $\ud$, but it is less obvious how to treat the codifferential operator $\ud^*$.

\subsubsection{Generalized Poisson problem}
Take for example the generalized Poisson problem using the Hodge-Laplacian acting on $k$-forms, $(\ud\ud^*+\ud^*\ud)u=f$, on $\Omega$ with boundary $\Gamma=\p\Omega$. A standard Galerkin approach, using integration by parts \eqref{integrationbyparts}, would give; find $u\in H\Lambda^k(\Omega)\cap H^*\Lambda^k(\Omega)$ with $\ud u\in H^*\Lambda^{k+1}(\Omega)$ and $\ud^*u\in H\Lambda^{k-1}(\Omega)$, given $f\in L^2\Lambda^k(\Omega)$, such that
	\begin{equation}
	\label{standardgalerkin}
	\big(\ud^* v,\ud^* u\big)_\Omega+\big(\ud v,\ud u\big)_\Omega=\big(v,f\big)_\Omega,\quad\quad\forall v\in H\Lambda^k(\Omega)\cap H^*\Lambda^k(\Omega).
	\end{equation}
It has a corresponding minimization problem for an energy functional over the space $H\Lambda^k(\Omega)\cap H^*\Lambda^k(\Omega)$. The standard Galerkin formulation is coercive, which immediately implies stability. Corresponding to this standard Galerkin formulation one usually chooses a $H^1\Lambda^k(\Omega)$-conforming approximation space. This could be a standard continuous piecewise polynomial vector space based on nodal interpolation.

However, in case of a nonconvex polyhedral or curvilinear or noncontractible domain $\Omega$, for allmost all $f$, $H\Lambda^k(\Omega)\cap H^*\Lambda^k(\Omega)\not\subset H^1\Lambda^k(\Omega)$. Consequently, the solution will be stable but inconsistent in general, \cite{costabel1991}. In other words, the solution converges to the wrong solution. Unfortunately, it seems not possible to construct $H\Lambda^k(\Omega)\cap H^*\Lambda^k(\Omega)$ conforming finite element spaces. Alternatively, one proposed to use \emph{mixed formulations}, \cite{brezzifortin}. In contrast to standard Galerkin, the mixed formulation uses integration by parts \eqref{integrationbyparts} to express each codifferential in terms of an exterior derivative and suitable boundary conditions.

Consequently, mixed formulations require only $H\Lambda^k(\Omega)$-conforming finite element spaces, which are much easier to construct. Therefore, in all cases mixed formulations do converge to the true solution. Mixed formulations correspond to saddle point problems instead of minimization problems.

The derivation of the mixed formulation of the Poisson problem consists of three steps: 
\begin{enumerate}
\item Introduce an auxiliary variable $\omega=\ud^*u$ in $H\Lambda^{k-1}$,
\item multiply both equations by test functions $\big(\tau,v\big)\in\{H\Lambda^{k-1}\times H\Lambda^{k}\}$ using $L^2$-inner products,
\item use integration by parts, as in \eqref{integrationbyparts}, to express the remaining codifferentials in terms of the exterior derivatives and boundary integrals.
\end{enumerate}
Again the boundary may constitute up to four different types of boundary conditions,
\begin{subequations}
\label{boundaryconditions2}
\begin{align}
\tr\omega&=0\quad\quad\  \mathrm{on}\ \Gamma_\omega,\\
\tr u&=0\quad\quad\  \mathrm{on}\ \Gamma_n,\\
\quad\tr\star u&=u_{b,t}\quad\ \mathrm{on}\ \Gamma_t,\quad\; \mathrm{with}\ u_{b,t}\in H^{\frac{1}{2}}\Lambda^{n-k}(\Gamma_t),\\
\tr\star\ud u&=g_b\quad\ \ \, \mathrm{on}\ \Gamma_\pi,\quad \mathrm{with}\ g_b\in H^{\frac{1}{2}}\Lambda^{n-k-1}(\Gamma_\pi).
\end{align}
\end{subequations}
Then also for the generalized Poisson problem the boundary $\Gamma$ consists up to four sections as defined in \eqref{boundarydecomposition}. To obtain a unique solution for the corresponding mized formulation, we define the following Hilbert spaces,
\begin{align}
W:=&\{\,\tau\in H\Lambda^{k-1}(\Omega)\,|\,\tr\tau=0\ \mathrm{on}\ \Gamma_\omega\,\},\\
V:=&\left\{
\begin{aligned}
&\{\, v\in H\Lambda^{k}(\Omega)\,|\,\tr v=0\ \mathrm{on}\ \Gamma_n\,\}\quad\mathrm{if}\ \Gamma_4=\emptyset,\\
&\{\, v\in H\Lambda^{k}(\Omega)\backslash\mathcal{Z}^{*,k}(\Omega)\,|\,\tr v=0\ \mathrm{on}\ \Gamma_n\,\}\quad\mathrm{if}\ \Gamma_4\neq\emptyset,
\end{aligned}
\right.
\end{align}
with corresponding norms, $\norm{\cdot}_W,\ \norm{\cdot}_V$, respectively. The resulting mixed formulation for the Poisson problem for all $0\leq k\leq n$ becomes: find $(\omega,u)\in\{W\times V\}$, given $f\in L^2\Lambda^k$, for all $(\tau,v)\in\{W\times V\}$, such that
\begin{subequations}
\label{poisson}
\begin{align}
\big(\tau,\omega\big)_\Omega-\big(\ud\tau,u\big)_\Omega&=-\int_{\Gamma_1\cup\Gamma_2}\tr\tau\wedge u_{b,t},\label{poisson1}\\ 
\big(v,\ud\omega\big)_\Omega+\big(\ud v,\ud u\big)_\Omega&=\big(v,f\big)_\Omega+\int_{\Gamma_2\cup\Gamma_4}\tr v\wedge g_b.\label{poisson2}  
\end{align}
\end{subequations}
Note that, for a scalar Poisson, it is not a choice whether to use Galerkin or mixed formulation, but it depends on whether the scalar is a 0-form or an $n$-form. This is determined by the physics.

\subsubsection{Stokes problem}
In a similar way the mixed formulation of the VVP formulation of the Stokes problem is obtained. Consider the Hilbert spaces $W$ and $V$ defined in the previous section, where $k=n-1$, and define the following Hilbert space
\begin{equation}
Q:=\left\{
\begin{aligned}
&q\in L^2\Lambda^n(\Omega),\quad\quad\quad\ \,\mathrm{if}\ \Gamma_\pi\neq\emptyset,\\
&q\in L^2\Lambda^n(\Omega)\backslash\mathcal{Z}^{*,n},\quad\mathrm{if}\ \Gamma_\pi=\emptyset,
\end{aligned}
\right.
\end{equation}
with corresponding norm $\norm{\cdot}_Q$ and where $\mathcal{Z}^{*,n}=\mathbb{R}$. Then the mixed formulation of the VVP formulation reads: find $(\omega,u,p)\in\{W\times V\times Q\}$, for the given data $f\in L^2\Lambda^{n-1}(\Omega)$, $g\in L^2\Lambda^{n}(\Omega)$ and natural boundary conditions $u_{b,t}\in H^{\frac{1}{2}}\Lambda^{n-1}(\Gamma_t)$, $\Pi_b\in H^{\frac{1}{2}}\Lambda^n(\Gamma_\pi)$, for all $\big(\tau,v,q\big)\in\{W\times V\times Q\}$, such that
\begin{subequations}
\label{mixedstokes}
\begin{align}
\big(\tau,\omega\big)_\Omega-\big(\ud\tau,u\big)_\Omega&=-\int_{\Gamma_1\cup\Gamma_2}\tr\tau\wedge u_{b,t},\label{mixedstokes1}\\
\big(v,\ud\omega\big)_\Omega+\big(\ud v,\ud u\big)_\Omega+\big(\ud v,p\big)_\Omega&=\big(v,f\big)_\Omega+\int_{\Gamma_2\cup\Gamma_4}\tr v\wedge\Pi_b,\label{mixedstokes2}\\
\big(q,\ud u\big)_\Omega&=\big(q,g\big)_\Omega.\label{mixedstokes3}
\end{align}
\end{subequations}
Again use is made of integration by parts, \eqref{integrationbyparts}. 
\begin{proposition}\cite{bernardi2006}
Problems \eqref{stokeseq}-\eqref{boundaryconditions} and \eqref{mixedstokes} are equivalent, in the sense that any triple $\big(\omega,u,p\big)\in\{W\times V\times Q\}$ is a solution of problem \eqref{stokeseq}-\eqref{boundaryconditions} if and only if it is a solution of problem \eqref{mixedstokes}.
\end{proposition}

\subsection{Well-posedness of mixed formulation}
Before we continue first define the following nullspaces of $W$,
\begin{subequations}
\begin{align}
Z_W&:=\{\;\tau\in W\;|\;\ud\tau=0\;\},\\
Z_W^*&:=\{\;\tau\in W\;|\;\ud^*\tau=0\;\},
\end{align}
and consider the following decompositions, $W=Z_W\oplus Z_W^\perp$ and $W=Z_W^*\oplus Z_W^{*,\perp}$. Since vorticity is defined as $\omega=\ud^*u$, we have $\omega\in Z_W^*$, and because we consider contractible domains only, it follows that $\omega\in Z_W^\perp$. Note that for $n=2$, $Z^*_W\equiv W$. A similar decomposition can be made for $V$. Define
\begin{equation}
Z_V:=\{\;v\in V\;|\;\ud v=0\;\},
\end{equation}
\end{subequations}
then $V=Z_V\oplus Z_V^\perp$. The velocity is decomposed as $u=u_\mathcal{Z}+u_\perp$, where $u_\mathcal{Z}\in Z_V$ and $u_\perp\in Z_V^\perp$.\\

We can write the mixed formulation of \eqref{mixedstokes} in a more general representation, using four continuous bilinear forms,
\begin{align*}
\sfa(\cdot,\cdot)&:=(\cdot,\cdot)_\Omega\ :W\times W\rightarrow\mathbb{R},\quad\; \sfb(\cdot,\cdot):=(\ud\cdot,\cdot)_\Omega\ :V\times Q\rightarrow\mathbb{R},\\ \sfc(\cdot,\cdot)&:=(\ud\cdot,\cdot)_\Omega\ :W\times V\rightarrow\mathbb{R},\quad \sfe(\cdot,\cdot):=(\ud\cdot,\ud\cdot)_\Omega\ :V\times V\rightarrow\mathbb{R},
\end{align*}
and three continuous linear forms
\begin{gather*}
\sff(\cdot):=(\cdot,f)_\Omega+\int_{\Gamma_2\cup\Gamma_4}\tr\cdot\wedge\Pi_b\ :V\rightarrow\mathbb{R},\\
\sfg(\cdot):=(\cdot,g)_\Omega\ :Q\rightarrow\mathbb{R},\quad\quad \sfh(\cdot):=-\int_{\Gamma_1\cup\Gamma_2}\tr\cdot\wedge u_{b,t}\ :W\rightarrow\mathbb{R}.
\end{gather*}
The mixed formulation becomes
\begin{subequations}
\label{doublesaddlepoint}
\begin{align}
\sfa(\tau,\omega)-\sfc(\tau,u)&=\sfh(\tau),\quad\quad\!\forall\tau\in W, \label{saddlepoint1}\\
\sfe(v,u)+\sfc(\omega,v)+\sfb(v,p)&=\sff(v),\quad\quad\forall v\in V, \label{saddlepoint2}\\
\sfb(u,q)&=\sfg(q).\quad\quad\forall q\in Q. \label{saddlepoint3}
\end{align}
\end{subequations}
There exists continuity constants $0<c_\sfa,c_\sfb,c_\sfc,c_\sfe<\infty$ such that
\begin{equation}
\sfa(\tau,\kappa)\leq c_\sfa\norm{\tau}_W\norm{\kappa}_W,\quad \sfb(v,q)\leq c_\sfb\norm{v}_V\norm{q}_Q,\quad \sfc(\tau,v)\leq c_\sfc\norm{\tau}_W\norm{v}_V,\quad\sfe(v,w)\leq c_\sfe\norm{v}_V\norm{w}_V.
\end{equation}
By Cauchy-Schwarz we know that $c_\sfa=1$, however we write $c_\sfa$ for generality purpose. The continuous linear forms are bounded such that
\begin{equation}
\sff(v)\leq \norm{f}\norm{v}_V,\quad \sfg(v)\leq \norm{g}\norm{v}_V,\quad\sfh(\tau)\leq\norm{h}\norm{\tau}_W.
\end{equation}
At first restrict to all $v=v_\mathcal{Z}\in Z_V$. This gives the vorticity-velocity subproblem, which is a saddle point problem:
\begin{subequations}
\label{vorticityvelocityproblem}
\begin{align}
\sfa(\tau,\omega)-\sfc(\tau,u_\mathcal{Z})&=\sfh(\tau),\quad\quad\ \,\forall\tau\in W,\\
\sfc(v_\mathcal{Z},\omega)&=\sff(v_\mathcal{Z}),\quad\quad\forall v_\mathcal{Z}\in Z_V.
\end{align}
\end{subequations}

\begin{proposition}\cite{dubois2002}\label{prop:infsupac}
System \eqref{vorticityvelocityproblem} has a unique solution $(\omega,u_\mathcal{Z})\in\{W\times Z_V\}$ if there exists positive constants $\alpha,\ \gamma$, such that
we have coercivity in the kernel of $W$,
\begin{equation}
\label{infsupa}
\inf_{\tau_\mathcal{Z}\in Z_W}\sup_{\kappa_\mathcal{Z}\in Z_W}\frac{\sfa(\tau_\mathcal{Z},\kappa_\mathcal{Z})}{\norm{\tau_\mathcal{Z}}_W\norm{\kappa_\mathcal{Z}}_W}\geq\alpha,\quad\quad
\inf_{\kappa_\mathcal{Z}\in Z_W}\sup_{\tau_\mathcal{Z}\in Z_W}\frac{\sfa(\tau_\mathcal{Z},\kappa_\mathcal{Z})}{\norm{\tau_\mathcal{Z}}_W\norm{\kappa_\mathcal{Z}}_W}\geq\alpha,
\end{equation}
and satisfies the following inf-sup condition for $\sfc(\tau,v_\mathcal{Z})$,

\begin{equation}
\label{infsupc}
\inf_{v_\mathcal{Z}\in Z_V}\sup_{\tau\in W}\frac{\sfc(\tau,v_\mathcal{Z})}{\norm{\tau}_W\norm{v_\mathcal{Z}}_V}\geq\gamma,
\end{equation}
\begin{proof}
The proof of \eqref{infsupa} is straightforward, see e.g. \cite{brezzifortin}. For \eqref{infsupc}, we have $\sfc(\tau,v_\mathcal{Z})=(\ud\tau,v_\mathcal{Z})_\Omega$, where $\ud:Z_W^\perp\rightarrow Z_V$. Thus, given $v_\mathcal{Z}\in Z_V$ there exists a unique $\tau_v\in Z_W^\perp$ such that $\ud\tau_v=v_\mathcal{Z}$ and $\norm{\tau_v}_W\leq c_P\norm{v_\mathcal{Z}}_V$ by \lemmaref{lem:poicareinequality}. Therefore
\[
\sup_{\tau\in W}\frac{\sfc(\tau,v_\mathcal{Z})}{\norm{\tau}_W}\geq\frac{\sfc(\tau_v,v_\mathcal{Z})}{\norm{\tau_v}_W}=\frac{\norm{v_\mathcal{Z}}_V^2}{\norm{\tau_v}_W}\geq\frac{1}{c_P}\norm{v_\mathcal{Z}}_V.
\]
\end{proof}
\end{proposition}
\begin{proposition}\label{prop:infsupeb}
The full problem \eqref{doublesaddlepoint} has a unique solution $(\omega,u,p)\in\{W,V,Q\}$ if conditions \eqref{infsupa} and \eqref{infsupc} from Proposition \ref{prop:infsupac} are satisfied and additionally is there exists positive constants $\ve,\beta$, such that we have coercivity in the range of $V$,
\begin{equation}
\label{infsupe}
\inf_{v_\perp\in Z_V^\perp}\sup_{w_\perp\in Z_V^\perp}\frac{\sfe(v_\perp,w_\perp)}{\norm{v_\perp}_V\norm{w_\perp}_V}\geq\ve,\quad\quad
\inf_{w_\perp\in Z_V^\perp}\sup_{v_\perp\in Z_V^\perp}\frac{\sfe(v_\perp,w_\perp)}{\norm{v_\perp}_V\norm{w_\perp}_V}\geq\ve,
\end{equation}
and satisfies the following inf-sup condition for $\sfb(v,q)$,
\begin{equation}
\label{infsupb}
\inf_{q\in Q}\sup_{v\in V}\frac{\sfb(v,q)}{\norm{v}_V\norm{q}_Q}\geq\beta>0.
\end{equation}
\begin{proof}
The proof is similar to that of \propref{prop:infsupac}. See also \cite{bochev2003}, Section 7.1.
\end{proof}
\end{proposition}
So well-posedness of the Stokes problem \eqref{mixedstokes} relies only on the Hodge decomposition and the Poincar\'e inequality.
\begin{corollary}\label{cor:wellposedness}\cite{bernardi2006,dubois2002}
Problem \eqref{mixedstokes} is well-posed according to Propositions \ref{prop:infsupac} and \ref{prop:infsupeb}. That is, for any given data $f\in L^2\Lambda^{n-1}(\Omega)$ and $g\in L^2\Lambda^n(\Omega)$ and natural boundary conditions $u_{b,t}\in H^{\frac{1}{2}}\Lambda^{n-1}(\Gamma_t)$ and $\Pi_b\in H^{\frac{1}{2}}\Lambda^n(\Gamma_\pi)$, there exists a unique solution $(\omega,u,p)\in W\times V\times Q$ satisfying \eqref{mixedstokes}. Moreover, this solution satisfies:
\begin{equation}
\norm{\omega}_{W}+\norm{u}_{V}+\norm{p}_{Q}\leq C\left(\norm{f}_{L^2\Lambda^{n-1}}+\norm{g}_{L^2\Lambda^n}+\norm{u_{b,t}}_{H^{\frac{1}{2}}\Lambda^{n-1}}+\norm{\Pi_b}_{H^{\frac{1}{2}}\Lambda^n}\right),
\end{equation}
where $C$ is a constant depending only on the Poincar\'e constant $c_P$ and the continuity constants.
\end{corollary}

\section{Compatible spectral discretization}\label{sec:discretization}
Well-posedness of the Stokes problem in VVP formulation relies solely on the Hodge decomposition and the Poincar\'e inequality. For a compatible discretization, these properties need to be respected as well in the finite dimensional spaces. Key ingredient to obtain a discrete Hodge decomposition and discrete Poincar\'e inequality is the construction of a bounded projection operator that commutes with the exterior derivative.

The compatible spectral discretization consists of three parts. First, the discrete structure is described in terms of chains and cochains from algebraic topology, the discrete counterpart of differential geometry. This discrete structure mimics many of the properties from differential geometry. Secondly, mimetic operators are introduced that relate the continuous formulation in terms of differential forms to the discrete representation based on cochains and finite dimensional differential forms. Thirdly, mimetic spectral element basis functions are described following the definitions of the mimetic operators. In this paper we address these topics only briefly. More details of the mimetic spectral element method can be found in \cite{kreeft2012,kreeftpalhagerritsma2011}. Finally, well-posedness of the discrete numerical formulation is proven and interpolation error estimates are given.

\subsection{Algebraic Topology}
Let $D$ be an \emph{oriented cell-complex} covering the manifold $\Omega$, describing the topology of the mesh, and consisting of $k$-cells $\tau_{(k)}$, $k=0,\hdots,n$. The two most popular classes of $k$-cells in literature to describe the topology of a manifold are either in terms of \emph{simplices}, see for instance \cite{munkres1984,singerthorpe,whitney}, or in terms of \emph{cubes}, see \cite{Massey2,tonti1}. From a topological point of view both descriptions are equivalent, see \cite{dieudonne}. Despite this equivalence of simplicial complexes and cubical complexes, the reconstruction maps in terms of basis functions, to be discussed in \secref{mimeticoperators}, differ significantly. For mimetic methods based on simplices see \cite{arnoldfalkwinther2006,arnoldfalkwinther2010,desbrun2005c,rapetti2009}, whereas for mimetic methods based on singular cubes see \cite{arnoldboffifalk2005,HymanShashkovSteinberg2002,hymansteinberg2004,RobidouxSteinberg2011}. We restrict ourselves to $k$-cubes, although we will keep calling them $k$-cells.

The ordered collection of all $k$-cells in $D$ generate a basis for the space of $k$-chains, $C_k(D)$. Then a $k$-chain, $\kchain{c}{k}\in C_k(D)$, is a formal linear combination of $k$-cells, $\tau_{(k),i}\in D$,
\begin{equation}
\label{kchain}
\kchain{c}{k}=\sum_ic_i\tau_{(k),i}.
\end{equation}
The boundary operator on $k$-chains, $\partial:\kchainspacedomain{k}{D}\spacemap\kchainspacedomain{k-1}{D}$, is an homomorphism defined by \cite{hatcher,munkres1984},
\begin{equation}\label{algebraic::boundary_operator}
\partial \kchain{c}{k} = \partial \sum_{i}c_{i}\tau_{(k),i} := \sum_{i}c_{i} \partial \left ( \tau_{(k),i} \right ) \;.
\end{equation}
The boundary of a $k$-cell $\tau_{(k)}$ will then be a $(k-1)$-chain formed by the oriented faces of $\tau_{(k)}$. Like the exterior derivative, applying the boundary operator twice on a $k$-chain gives the null $(k-2)$-chain, $\p\p\kchain{c}{k}=\kchain{0}{k-2}$ for all $\kchain{c}{k}\in C_k(D)$. The set of $k$-chains and boundary operators gives rise to an exact sequence, the chain complex $(C_k(D),\p)$,
\begin{equation}
\label{chaincomplex}
\begin{CD}
\cdots @<\p<< C_{k-1}(D) @<\p<< C_k(D) @<\p<< C_{k+1}(D) @<\p<< \cdots.
\end{CD}
\end{equation}
Let $B_k$ be the range and $Z_k$ be the nullspace of $\p$ in $C_k$. Then the topological Hodge decomposition of the space of $k$-chains is given by $C_k=Z_k\oplus Z_k^\perp$, where $Z_k=B_k$ on contractible domains\footnote{Although `perpendicular' in a topological space is not well defined, we refer to $Z^\perp_k$ as the complement space of $Z_k$ in $C_k$.}. The boundary operator on chains in \eqref{chaincomplex} is a bijection that maps $\p:Z_k^\perp\rightarrow B_{k-1}$.

Dual to the space of $k$-chains, $C_k(D)$, is the space of \emph{$k$-cochains}, $C^k(D)$, defined as the set of all linear functionals, $\kcochain{c}{k}:C_k(D)\rightarrow\mathbb{R}$. The duality is expressed using the duality pairing $\langle\kcochain{c}{k},\kchain{c}{k}\rangle:=\kcochain{c}{k}(\kchain{c}{k})$. Note the resemblance between this duality pairing and the integration of differential forms \eqref{integration}.

Let $\{\tau_{(k),i}\}$ form a basis of $C_k(D)$, then there is a dual basis $\{\tau^{(k),i}\}$ of $C^k(D)$, such that $\tau^{(k),i}(\tau_{(k),i})=\delta^i_j$ and all $k$-cochains can be represented as linear combinations of the basis elements,
\begin{equation}
\kcochain{c}{k}=\sum_ic^i\tau^{(k),i}.
\end{equation}
With the duality relation between chains and cochains, we can define the formal adjoint of the boundary operator which constitutes an exact sequence on the spaces of $k$-cochains in the cell complex. This formal adjoint is called the \emph{coboundary operator}, $\delta:\kcochainspacedomain{k}{D}\rightarrow\kcochainspacedomain{k+1}{D}$, and is defined analogous to \eqref{stokestheorem} as
\begin{equation}
\duality{\delta\kcochain{c}{k}}{\kchain{c}{k+1}} := \duality{\kcochain{c}{k}}{\partial\kchain{c}{k+1}}, \quad\forall\kcochain{c}{k}\in\kcochainspacedomain{k}{D} \text{ and  } \,\forall\kchain{c}{k+1}\in\kchainspacedomain{k+1}{D} \;. \label{coboundary}
\end{equation}
Note that expression \eqref{coboundary} is nothing but a discrete Stokes' theorem and that the coboundary operator is nothing but a discrete exterior derivative. Also the coboundary operator satisfies $\delta\delta\kcochain{c}{k}=\kcochain{0}{k+2}$, for all $\kcochain{c}{k}\in C^k(D)$, and gives rise to an exact sequence, called the \emph{cochain complex} $(C^k(D),\delta)$,
\begin{equation}
\label{cochaincomplex}
\begin{CD}
\cdots @>\delta>> C^{k-1}(D)@>\delta>> C^k(D)@>\delta>>C^{k+1}(D)@>\delta>>\cdots\;.
\end{CD}
\end{equation}
Let $B^k$ be the range and $Z^k$ be the nullspace of $\delta$ in $C^k$, then a Hodge decomposition of the space of $k$-cochains is given by $C^k=Z^k\oplus Z^{k,\perp}$, where $Z^k=B^k$ on contractible domains. The coboundary operator in \eqref{cochaincomplex} is a bijection that maps $\delta:Z^{k,\perp}\rightarrow B^{k+1}$. Note the similarity between this map, that of the boundary operator on $k$-chains and that of the exterior derivative on $k$-forms.

\subsection{Mimetic Operators} \label{mimeticoperators}
The discretization of the flow variables involves a bounded projection operator, $\pi_h$, from the complete space $H\Lambda^k(\Omega)$ to a conforming subspace $\Lambda^k_h(\Omega;C_k)\subset H\Lambda^k(\Omega)$.
The projection operation consists of two steps, a \emph{reduction operator}, $\reduction:H\Lambda^k(\Omega)\rightarrow C^k(D)$, that integrates the $k$-forms on $k$-chains to get $k$-cochains, and a \emph{reconstruction operator}, $\reconstruction:C^k(D)\rightarrow\Lambda^k_h(\Omega;C_k)$, to reconstruct $k$-forms from $k$-cochains using appropriate basis-functions. These mimetic operators were already introduced before in \cite{bochevhyman2006,hymanscovel1988}. A composition of the two gives the projection operator $\projection=\reconstruction\circ\reduction$ as is illustrated below.
\begin{diagram}
H\Lambda^k(\Omega)&\rTo^{\quad\projection\quad}& \Lambda^k_h(\Omega;C_k)\\
\dTo^{\reduction} & \ruTo_{\reconstruction} & \\
C^k(D) & &
\end{diagram}
These three operators together constitute the mimetic framework. An extensive discussion on mimetic operators can be found in \cite{kreeft2012,kreeftpalhagerritsma2011}.

The reduction $\reduction$ and reconstruction $\reconstruction$ operators are defined below. The fundamental property of $\reduction$ and $\reconstruction$ is the commutation with differentiation in terms of exterior derivative and coboundary operator.

The reduction operator $\reduction:H\kformspacedomain{k}{\Omega}\rightarrow \kcochainspacedomain{k}{D}$ is a homomorphism that maps differential forms to cochains. This map is defined by integration as
\begin{equation}
\duality{\reduction a}{\tau_{(k)}}:=\int_{\tau_{(k)}}a,\quad \forall a\in H\Lambda^k(\Omega),\ \tau_{(k)}\in C_k(D).
\label{reduction}
\end{equation}
Then for all $\kchain{c}{k}\in C_k(D)$, the reduction of the $k$-form, $a\in H\Lambda^{k}(\Omega)$, to the $k$-cochain, $\kcochain{a}{k}\in C^k(D)$, is given by
\begin{equation}
\kcochain{a}{k}(\kchain{c}{k}):=\duality{\reduction a}{\kchain{c}{k}}\stackrel{\eqref{kchain}}{=}\sum_ic^i\duality{\reduction a}{\tau_{(k),i}}\stackrel{\eqref{reduction}}{=}\sum_ic^i\int_{\tau_{(k),i}}a=\int_{\kchain{c}{k}}a.
\end{equation}

The reduction maps has a commuting property with respect to differentiation in terms of exterior derivative and coboundary operator,
\begin{equation}
\reduction\ud=\delta\reduction,\quad\mathrm{on}\ H\Lambda^k(\Omega).\label{cdp1}
\end{equation}
Since $\reduction$ is defined by integration, \eqref{cdp1} follows directly from Stokes theorem \eqref{stokestheorem} and the duality property \eqref{coboundary}.

Next by definition also the reconstruction map $\reconstruction:C^k(D)\rightarrow\Lambda^k_h(\Omega;C_k)$ needs to have a commuting property with respect to differentiation in terms of exterior derivative and coboundary operator,
\begin{equation}
\ud\reconstruction=\reconstruction\delta,\quad\mathrm{on}\ C^k(D).\label{cdp2}
\end{equation}
The reconstruction $\reconstruction$ must be the right inverse of $\reduction$, so $\mathcal{RI}=Id$ on $C^k(D)$, and we want it to be an approximate left inverse of $\mathcal{R}$, so $\mathcal{IR}=Id+\mathcal{O}(h^p)$ on $H\Lambda^k(\Omega)$. This composition is defined as the projection operator.
\begin{definition}[\textbf{Bounded projection operator}]\label{def:projection}
The composition $\reconstruction\circ\reduction$ will denote the projection operator, $\projection\define\reconstruction\reduction:H\kformspace{k}(\Omega)\rightarrow\kformspace{k}_h(\Omega;C_k)$, allowing for a finite dimensional representation of a $k$-form,
\begin{equation}
\projection a:=\reconstruction\reduction a, \quad \pi_ha\in\kformspace{k}_h(\Omega;C_k)\subset H\kformspace{k}(\Omega).
\label{projection}
\end{equation}
where $\reconstruction\reduction a$ is expressed as a combination of $k$-cochains and interpolating $k$-forms. The projection operator $\pi_h$ is a bounded operator if for $C<\infty$ and for all $a\in\Lambda^k(\Omega)$ we have $\norm{\pi_ha}_{H\Lambda^k}\leq C\norm{a}_{H\Lambda^k}$.
\end{definition}
A proof that $\pi_h$ is indeed a projection operator is given in \cite{kreeftpalhagerritsma2011}. In \secref{sec:boundedprojection} also boundedness is proven.
\begin{lemma}[\textbf{Commutation property}]
\label{Lem:projectionextder}
There exists a commuting property for the projection and the exterior derivative, such that
\begin{equation}
\label{projectionextder}
\ederiv\projection=\projection\ederiv\quad\mathrm{on}\ H\Lambda^k(\Omega).
\end{equation}
\begin{proof}
Express the projection in terms of the reduction and reconstruction operator, then
\[
\ederiv\projection\difform{a}\stackrel{\eqref{projection}}{=}\ederiv\reconstruction\reduction\difform{a}\stackrel{\eqref{cdp2}}{=}\reconstruction\delta\reduction\difform{a}\stackrel{\eqref{cdp1}}{=}\reconstruction\reduction\ederiv\difform{a}\stackrel{\eqref{projection}}{=}\projection\ederiv\difform{a},\quad\forall a\in H\Lambda^k(\Omega).
\]
\end{proof}
\end{lemma}
Note that it is the intermediate step $\mathcal{I}\delta\mathcal{R}a$ that is used in practice for the discretization, see \cite{kreeft2012}, and \eqref{edgediscretization} on page \pageref{edgediscretization}.
\begin{corollary}[\textbf{Discrete Hodge decomposition}]\label{cor:projectionextder}
From \lemmaref{Lem:projectionextder} it follows that $\mathcal{B}^k_h:=\pi_h\mathcal{B}^k\subset\mathcal{B}^k$, $\mathcal{Z}^k_h:=\pi_h\mathcal{Z}^k\subset\mathcal{Z}^k$ and that on contractible domains, $\mathcal{Z}^{k,\perp}_h:=\pi_h\mathcal{Z}^{k,\perp}\subset\mathcal{Z}^{k,\perp}$. Then the discrete Hodge decomposition is given by $\Lambda^k=\mathcal{Z}^k_h\oplus\mathcal{Z}^{k,\perp}_h$. As a consequence of \lemmaref{Lem:projectionextder} and the discrete Hodge decomposition we have $Z_{W_h}\subset Z_W$, $Z_{V_h}\subset Z_V$, $\ud W_h \subset V_h$ and $\ud V_h= Q_h$, which shows that the discretization method is \emph{compatible}.
\end{corollary}

Finally, we do not restrict ourselves to affine mappings only, as is required in many other compatible finite elements, like N\'ed\'elec and Raviart-Thomas elements and their generalizations \cite{arnoldfalkwinther2006,nedelec1980,raviartthomas1977}, but also allow non-affine maps such as curvilinear transfinite or isoparametric mappings of quadrilaterals or hexahedrals, \cite{gordonhall1973}, where $\Phi$ and its inverse are piecewise sufficiently smooth, i.e.
\begin{enumerate}
\item $\Phi$ is a $\mathcal{C}^{p+1}$-diffeomorphism,
\item $|\Phi|_{W_\infty^l}\leq Ch^l,\quad\quad l\leq p+1$,
\item $|\Phi^{-1}|_{W_\infty^l}\leq Ch^{-l},\quad l\leq p+1$.
\end{enumerate}
This allows for better approximations in complex domains with curved boundaries, without the need for excessive refinement, while maintaining design convergence rates, \cite{ciarletraviart1972}. This is possible since the projection operator $\pi_h$ commutes with the pullback $\pullback$,
\begin{equation}
\label{pullbackprojection}
\pullback\projection=\projection\pullback\quad\mathrm{on}\ H\Lambda^k(\Omega).
\end{equation}
An extensive proof is given in \cite{kreeftpalhagerritsma2011}.

\subsection{Numerical stability}
Essential ingredients in proving numerical stability are the discrete Hodge decomposition and the discrete Poincar\'e inequality. Because the complexes $(H\Lambda,\ud)$ and $(\Lambda_h,\ud)$ are each others supercomplex and subcomplex, respectively, the discrete Poincar\'e inequality is directly related to the Poincar\'e inequality in \lemmaref{lem:poicareinequality} and the bounded projection in \defref{def:projection}.
\begin{lemma}[\textbf{Discrete Poincar\'e inequality}]\label{lem:discretepoincareinequality}Let $(H\Lambda,\ud)$ be a bounded closed Hilbert complex, $(\Lambda_h,\ud)$ a subcomplex, and $\pi_h$ a bounded projection. Then
\begin{equation}
\left\{
\begin{aligned}
&\norm{a_h}_{H\Lambda^k}\leq c_{Ph}\norm{\ud a_h}_{L^2\Lambda^k},\quad a_h\in\mathcal{Z}^{k,\perp}_h,\\
&1\leq c_{Ph}\leq c_P.
\end{aligned}
\right.
\end{equation}
\begin{proof} Given $a_h\in\mathcal{Z}_h^{k,\perp}$. From the Hodge decomposition and the bounded projection it follows that $\mathcal{B}_h^k\subset\mathcal{B}^k$ and $\mathcal{Z}^k_h\subset\mathcal{Z}^k$ and from the commutation relation \eqref{projectionextder} it follows that $\mathcal{Z}^{k,\perp}_h\subset\mathcal{Z}^{k,\perp}$. Since we consider a proper subspace, \lemmaref{lem:poicareinequality} is still valid, with $c_{Ph}\leq c_P$.
\end{proof}
\end{lemma}

\begin{theorem}[\textbf{Discrete well-posedness}]\label{th:discretewellposedness}
Let $(\Lambda_h,\ud)$ be a subcomplex of the closed Hilbert complex $(H\Lambda,\ud)$. Then there exists constants $\alpha_h,\beta_h,\gamma_h$, depending only on $c_{Ph}$, such that for any $(\tau_h,v_h,q_h)\in W_h\times V_h\times Q_h$, there exists a stable finite dimensional solution $(\omega_h,u_h,p_h)\in W_h\times V_h\times Q_h$ of the Stokes problem \eqref{mixedstokes}, with
\begin{equation}
\alpha_h>\alpha>0,\quad\beta_h>\beta>0,\quad\gamma_h>\gamma>0.
\end{equation}
\begin{proof}
This is just Propositions \ref{prop:infsupac} and \ref{prop:infsupeb} applied to the complex $(\Lambda_h,\ud)$, combined with the fact that the constant in the Poincar\'e inequality for $\Lambda^k_h$ is $c_{Ph}\leq c_P$ by \lemmaref{lem:discretepoincareinequality}.
\end{proof}
\end{theorem}

\subsection{Mimetic spectral element basis-functions}\label{sec:mimeticsem}
The finite dimensional differential forms used in this paper are polynomials, based on the idea of spectral element methods, \cite{canuto1}. The mimetic spectral elements used here were derived independently in \cite{gerritsma2011,robidoux2008}, and are more extensively discussed in \cite{kreeftpalhagerritsma2011}. Only the most important properties of the mimetic spectral element method are presented here.

In spectral element methods the domain $\Omega$ is decomposed into $M$ non-overlapping, in this case curvilinear quadrilateral or hexahedral, closed sub-domains $Q_m$,
\[
\Omega=\bigcup_{m=1}^MQ_m,\quad Q_m\cap Q_l=\p Q_m\cap\p Q_l,\ m\neq l,
\]
where in each sub-domain a Gauss-Lobatto grid is constructed. The complete mesh is indicated by $\mathcal{Q}:=\sum_{m=1}^MQ_m$.

The collection of Gauss-Lobatto meshes in all elements $Q_m\in\mathcal{Q}$ constitutes the cell complex $D$. For each element $Q_m$ there exists a sub cell complex, $D_m$. Note that $D_m\cap D_l,\ m\neq l$, is not an empty set in case they are neighboring elements, but contains all $k$-cells, $k<n$, of the common boundary.

Each sub-domain $Q_m\in\mathcal{Q}$ is mapped from the reference element, $\widehat{Q}=[-1,1]^n$, using the mapping $\Phi_m:\widehat{Q}\rightarrow Q_m$. Then all flow variables defined on $Q_m$ are pulled back onto this reference element using the following pullback operation, $\Phi^\star_m:\Lambda^k_h(Q_m;C_k)\rightarrow\Lambda^k_h(\widehat{Q};C_k)$.

The basis-functions that interpolate the cochains on the quadrilateral or hexahedral elements are constructed using tensor products. It is therefore sufficient to derive interpolation functions in one dimension and use tensor products afterwards to construct $n$-dimensional basis functions. A similar approach was taken in \cite{buffa2011b}. Because projection operator and pullback operator commute \eqref{pullbackprojection}, the interpolation functions are discussed for the reference element only. Since the mappings $\Phi_m$ and their inverse are assumed to be sufficiently smooth, the rates of convergence for interpolation estimates on the physical elements are equal to that of the reference element. Only the constants $C$ that will appear below will depend on the mappings $\Phi_m$, but will be independent of the meshsize and polynomial order.\\

\begin{figure}[htb]
      \begin{minipage}[t]{0.49\linewidth}
            \centering\includegraphics[width=0.9\linewidth]{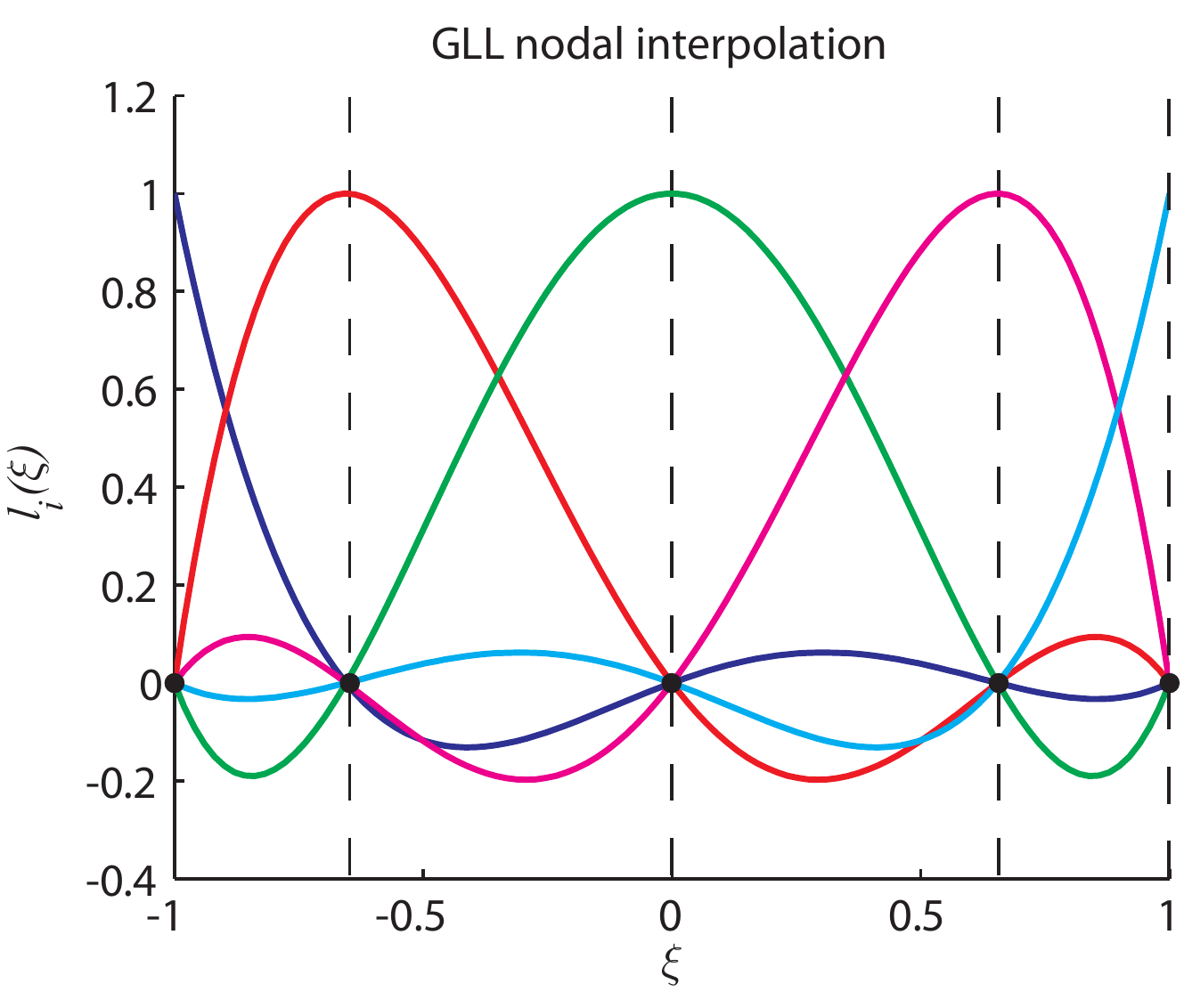}
	\caption{Lagrange polynomials on Gauss-Lobatto-Legendre grid.}
	\label{fig:lagrange}
    \end{minipage}\hfill
    \begin{minipage}[t]{0.49\linewidth}
            \centering\includegraphics[width=0.9\linewidth]{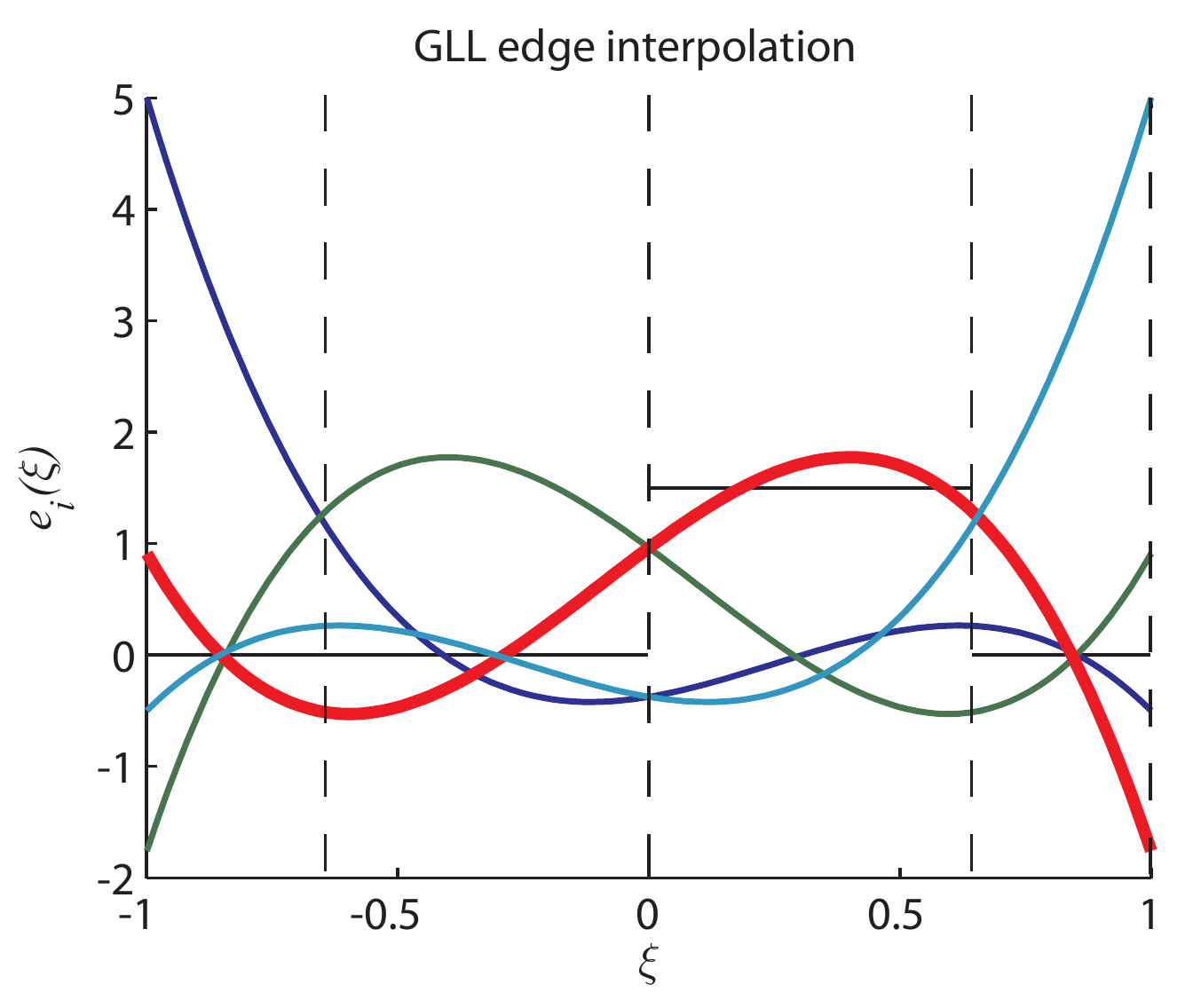}
	\caption{Edge polynomials on Gauss-Lobatto-Legendre grid.}
	\label{fig:edgepoly}
    \end{minipage}
\end{figure}

Consider a 0-form $a\in H\Lambda^0(\widehat{Q})$ on $\widehat{Q}:=\xi\in[-1,1]$, on which a cell complex $D$ is defined that consists of $N+1$ nodes, $\xi_i$, where $-1\leq \xi_0<\hdots<\xi_N\leq 1$, and $N$ edges, $\tau_{(1),i}=[\xi_{i-1},\xi_i]$, of which the nodes are their boundaries. Corresponding to this set of nodes (0-chains) there exists a projection using $N^{\rm th}$ order \textit{Lagrange polynomials}, $l_i(\xi)$, to approximate a $0$-form, as
\begin{equation}
\projection a=\sum_{i=0}^N a_il_i(\xi).
\label{nodalapprox}
\end{equation}
The property of Lagrange polynomials is that they interpolate nodal values. They are therefore suitable to reconstruct a 0-form form the 0-cochain $\kcochain{a}{0}=\reduction a$, $a\in\Lambda^0(\Omega)$, containing the set $a_i=a(\xi_i)$ for $i=0,\hdots,N$. Lagrange polynomials are in fact 0-forms, $l_i(\xi)\in\Lambda^0_h(\widehat{Q};C_0)$. Lagrange polynomials are constructed such that their value is one in the corresponding point and zero in all other grid points,
\begin{equation}
 \reduction l_i(\xi)=l_i(\xi_p)=\left\{
\begin{aligned}
&1&{\rm if}\ i=p\\ &0&{\rm if}\ i\neq p
\end{aligned}
\right..
 \label{nodalproperty}
\end{equation}
In \cite{gerritsma2011,robidoux2008} a similar basis for projection of 1-forms was derived, consisting of $1$-cochains and $1$-form polynomials, that is called the \textit{edge polynomial}, $e_i(\xi)\in\Lambda^1_h(\widehat{Q})$. Let $b\in L^2\Lambda^1(\widehat{Q})$, then the projected 1-form is given by
\begin{equation}
\projection b(\xi)=\sum_{i=1}^Nb_i e_i(\xi),
\end{equation}
where the edge polynomial is defined as
\begin{align}
\label{edge}
e_i(\xi)&=-\sum_{k=0}^{i-1}\ederiv l_k(\xi)=\sum_{k=i}^{N}\ederiv l_k(\xi)=\tfrac{1}{2}\sum_{k=i}^{N}\ederiv l_k(\xi)-\tfrac{1}{2}\sum_{k=0}^{i-1}\ederiv l_k(\xi).
\end{align}
Let $a_h(\xi)\in\Lambda^0_h(\widehat{Q};C_0)$ be expressed as in \eqref{nodalapprox}, then $b_h=\ud a_h\in\Lambda^1_h(\widehat{Q};C_1)$ is expressed as
\begin{equation}
\label{edgediscretization}
\ud\sum_{i=0}^Na_il_i(\xi)=\sum_{i=1}^N(a_i-a_{i-1})e_i(\xi)=\sum_{i=1}^N\big(\delta\kcochain{a}{0}\big)_ie_i(\xi)=\sum_{i=1}^Nb_ie_i(\xi),
\end{equation}
where $\delta$ is the coboundary operator \eqref{coboundary}, applied to the 0-cochain $\kcochain{a}{0}$.
It therefore satisfies \eqref{cdp2}. For derivations and proofs see \cite{gerritsma2011,kreeftpalhagerritsma2011,robidoux2008}. Similar to \eqref{nodalproperty}, the edge basis-functions are constructed such that when integrating $e_i(\xi)$ over a line segment it gives one for the corresponding element and zero for any other line segment, so
\begin{equation}
 \reduction e_i(\xi)=\int_{\xi_{p-1}}^{\xi_p}e_i(\xi)=\left\{
\begin{aligned}
&1&{\rm if}\ i=p\\ &0&{\rm if}\ i\neq p
\end{aligned}
\right..
 \label{intedge}
\end{equation}
Equations \eqref{nodalproperty} and \eqref{intedge} show that indeed we have $\mathcal{RI}=Id$. The fourth-order Lagrange and third-order edge polynomials, corresponding to a Gauss-Lobatto grid with $N=4$, are shown in Figures \ref{fig:lagrange} and \ref{fig:edgepoly}.

\subsection{Bounded projections and interpolation estimates}\label{sec:boundedprojection}
The mimetic framework uses Lagrange, $l_i(\xi)\in H\Lambda^0(\widehat{Q})$, and edge functions, $e_i(\xi)\in L^2\Lambda^1(\widehat{Q})$, for the reconstruction, $\mathcal{I}$.
Because we consider tensor products to construct higher-dimensional interpolation, it is sufficient to show that the projection operator is bounded in one dimension. A similar approach was used in \cite{buffa2011b}. Due to the way the edge functions are constructed, there exists a commuting diagram property between projection and exterior derivative,
\[
\begin{CD}
\mathbb{R} @>>> H\Lambda^0 @>\ederiv>> L^2\Lambda^1 @>>>0 \\
@.  @VV \pi_h V @VV\pi_hV  @. \\
\mathbb{R} @>>> \Lambda_h^0 @>\ederiv>> \Lambda_h^1 @>>>0,
\end{CD}
\]
which gives, for $a\in H\Lambda^0(\widehat{Q})$, the one form $\ederiv\pi_ha=\pi_h\ederiv a$ in $L^2\Lambda^1(\widehat{Q})$.
Lagrange interpolation by itself does not guarantee a convergent approximation \cite{erdos1980}, but it requires a suitably chosen set of points, $-1\leq\xi_0<\xi_1<\hdots<\xi_N\leq1$. Here, the Gauss-Lobatto distribution is proposed, because of its superior convergence behaviour. For $a\in H^m\Lambda^0(\Omega)$, the a priori error estimate in the $H\Lambda^0$-norm is given by \cite{canuto1},
\begin{equation}
\label{hpconvergenceestimate}
\Vert a-\pi_ha\Vert_{H\Lambda^0}\leq Ch^l|a|_{H^{l+1}\Lambda^0},\quad l=\mathrm{min}(N,m-1).
\end{equation}
Equation \eqref{hpconvergenceestimate} also implies that the projection of zero-forms is stable in the $H\Lambda^0(\widehat{Q})$, as is shown in the following proposition.

\begin{proposition}\label{boundedh}\cite{kreeftpalhagerritsma2011}
For $a\in H\Lambda^0(\widehat{Q})$ and the projection $\pi_h:H\Lambda^0\rightarrow \Lambda^0_h$, there exists the following two stability estimates in $H\Lambda^0$-norm and $H\Lambda^0$-semi-norm:
\begin{align}
\Vert\pi_ha\Vert_{H\Lambda^0}&\leq C\Vert a\Vert_{H\Lambda^0},\label{H1normstability}\\
|\pi_ha|_{H\Lambda^0}&\leq C|a|_{H\Lambda^0}.\label{H1seminormstability}
\end{align}
\end{proposition}
Now that we have a bounded linear projection of zero forms in one dimension, we can also proof boundedness of the projection of one-forms.
\begin{proposition}\label{boundede}
Let $a\in H\Lambda^0$ and $b=\ederiv a\in L^2\Lambda^1$, then there exists a bounded linear projection $\pi_h:L^2\Lambda^1\rightarrow\Lambda^1_h$, such that
\begin{equation}
\Vert\pi_hb\Vert_{L^2\Lambda^1}\leq C\Vert b\Vert_{L^2\Lambda^1}.
\end{equation}
\begin{proof}
The proof is based on the result of the previous proposition and the commutation between the bounded projection operator and the exterior derivative, \lemmaref{Lem:projectionextder},
\[
\Vert\pi_hb\Vert_{L^2\Lambda^1}=|\pi_h\ederiv a|_{L^2\Lambda^1}=|\ederiv\pi_ha|_{L^2\Lambda^1}=|\pi_ha|_{H\Lambda^0}\leq C|a|_{H\Lambda^0}=C\Vert\ederiv a\Vert_{L^2\Lambda^1}=C\Vert b\Vert_{L^2\Lambda^1}.
\]
\end{proof}
\end{proposition}
Propositions \ref{boundedh} and \ref{boundede} show that the projection $\pi_h$ is a \emph{bounded projection operator}, based on Lagrange functions and edge functions. As for zero forms using Lagrange interpolation, we can also give an estimate for the interpolation error of one forms, interpolated using edge functions.
\begin{proposition}\cite{kreeftpalhagerritsma2011}\label{prop:convratee}
Let $a\in H\Lambda^0$ and $b=\ederiv a\in L^2\Lambda^1$, the interpolation error $b-\pi_hb\in L^2\Lambda^1$ is given by
\begin{equation}
\label{edgeinterpolationerror}
\Vert b-\pi_hb\Vert_{L^2\Lambda^1}\leq Ch^{l}|b|_{H^{l}\Lambda^1},\quad l=\mathrm{min}(N,m-1).
\end{equation}
\end{proposition}

The one dimensional results can be extended to the multidimensional framework by means of tensor products. This allows for the interpolation of integral quantities defined on $k$-dimensional cubes. Consider a reference element in $\mathbb{R}^2$, $\widehat{Q}=[-1,1]^2$. Then the interpolation functions for points, lines, surfaces (2D volumes) are given by,
\begin{displaymath}
\begin{aligned}
&\mathrm{point}:&&P^{(0)}_{i,j}(\xi,\eta)=l_i(\xi)\otimes l_j(\eta),\\
&\mathrm{line}:&&L^{(1)}_{i,j}(\xi,\eta)=\{e_i(\xi)\otimes l_j(\eta),\ l_i(\xi)\otimes e_j(\eta)\},\\
&\mathrm{surface}:&&S^{(2)}_{i,j}(\xi,\eta)=e_i(\xi)\otimes e_j(\eta).
\end{aligned}
\end{displaymath}
The approximation spaces are spanned by combinations of Lagrange and edge basis functions,
\begin{align*}
H^1\Lambda^0(\Omega)\supset\Lambda^0_h(\mathcal{Q};C_0)&:=\mathrm{span}\left\{P^{(0)}_{i,j}\right\}_{i=0,j=0}^{N,N},\\
H\Lambda^1(\Omega)\supset\Lambda^1_h(\mathcal{Q};C_1)&:=\mathrm{span}\left\{\big(L^{(1)}_{i,j}\big)_1\right\}_{i=1,j=0}^{N,N}\times \mathrm{span}\left\{\big(L^{(1)}_{i,j}\big)_2\right\}_{i=0,j=1}^{N,N},\\
L^2\Lambda^2(\Omega)\supset\Lambda^2_h(\mathcal{Q};C_2)&:=\mathrm{span}\left\{S^{(2)}_{i,j}\right\}_{i=1,j=1}^{N,N}.
\end{align*}
For the variables vorticity, velocity and pressure in the VVP formulation of the Stokes problem, the $h$-convergence rates of the interpolation errors in $L^2\Lambda^k$-norm become,
\begin{equation}
\norm{\omega-\pi_h\omega}_{L^2\Lambda^{n-2}}=\mathcal{O}(h^{N+s}),\quad\norm{u-\pi_hu}_{L^2\Lambda^{n-1}}=\mathcal{O}(h^N),\quad\norm{p-\pi_hp}_{L^2\Lambda^n}=\mathcal{O}(h^N),
\label{L2interpolation}
\end{equation}
in case the functions $(\omega,u,p)$ are sufficiently smooth, where $s=1$ for $n=2$ and $s=0$ for $n>2$. The interpolation errors in $H\Lambda^k$-norm become,
\begin{equation}
\norm{\omega-\pi_h\omega}_{H\Lambda^{n-2}}=\mathcal{O}(h^N),\quad\norm{u-\pi_hu}_{H\Lambda^{n-1}}=\mathcal{O}(h^N),
\label{Hinterpolation}
\end{equation}
with $N$ defined as in \secref{sec:mimeticsem}.

\section{Error estimates}\label{sec:errorestimates}
Next consider the finite dimensional problem: find $(\omega_h,u_h,p_h)\in\{W_h\times V_h\times Q_h\}$, given $f\in L^2\Lambda^{n-1}(\Omega)$ and $g\in L^2\Lambda^{n}$ and boundary conditions in \eqref{boundaryconditions}, for all $(\tau_h,v_h,q_h)\in\{W_h\times V_h\times Q_h\}$, such that
\begin{subequations}
\label{discretedouble}
\begin{align}
\sfa(\tau_h,\omega_h)-\sfc(\tau_h,u_h)&=\sfh(\tau_h),\quad\quad\forall\tau_h\in W_h, \label{discretedouble1}\\
\sfe(v_h,u_h)+\sfc(\omega_h,v_h)+\sfb(v_h,p_h)&=\sff(v_h),\quad\quad\forall v_h\in V_h, \label{discretedouble2}\\
\sfb(u_h,q_h)&=\sfg(q_h).\quad\quad\forall q_h\in Q_h. \label{discretedouble3}
\end{align}
\end{subequations}
The following theorem gives the a priori error estimates of this problem when using the compatible spectral discretization method described in the previous section. Corollary~\ref{cor:projectionextder} showed that we have $Z_{W_h}\subset Z_W$ and $Z_{V_h}\subset Z_V$. From this it follows that we have compatible finite dimensional subspaces: $W_h\subset W$, $V_h=\ud W_h\oplus\ud^*Q_h\subset V$ and $Q_h=\ud V_h\subset Q$. The derivations of the error estimates are based on the methodology of \cite{brezzifortin}. The proofs are given in the subsequent propositions.
\begin{theorem}[\textbf{Error estimates}]\label{th:errorestimates}
Let $(\omega,u,p)$ be the solution of the continuous problem given in \eqref{mixedstokes} or \eqref{doublesaddlepoint} and $(\omega_h,u_h,p_h)$ the solution of the finite dimensional problem in \eqref{discretedouble}. The continuous problem is well-posed by Propositions \ref{prop:infsupac} and \ref{prop:infsupeb} and the finite dimensional problem is well-posed by Theorem~\ref{th:discretewellposedness} and Propositions~\ref{boundedh} and \ref{boundede}. Furthermore, from Corollary~\ref{cor:projectionextder} we have that for the compatible spectral discretization method, $Z_{W_h}\subset Z_W$ and $Z_{V_h}\subset Z_V$. Then the following a priori error estimates for the VVP formulation of the Stokes problem hold:
\begin{gather}
\norm{\omega-\omega_h}_W\leq\left(1+\frac{c_\sfa}{\alpha_h}\right)\left(1+\frac{c_\sfc}{\gamma_h}\right)\inf_{\tau_h\in W_h}\norm{\omega-\tau_h}_W,\\
\norm{u-u_h}_V\leq\left(1+\frac{c_\sfc}{\gamma_h}\right)\left(1+\frac{c_\sfb}{\beta_h}\right)\inf_{v_h\in V_h}\norm{u-v_h}_V+\frac{c_\sfa}{\gamma_h}\left(1+\frac{c_\sfa}{\alpha_h}\right)\left(1+\frac{c_\sfc}{\gamma_h}\right)\inf_{\tau_h\in W_h}\norm{\omega-\tau_h}_W,\\
\norm{p-p_h}_Q\leq\left(1+\frac{c_\sfb}{\beta_h}\right)\inf_{q_h\in Q_h}\norm{p-q_h}_Q+\frac{c_\sfe}{\beta_h}\left(1+\frac{c_\sfc}{\gamma_h}\right)\left(1+\frac{c_\sfb}{\beta_h}\right)\inf_{v_h\in V_h}\norm{u-v_h}_V\\
+\left(\frac{c_\sfa}{\gamma_h}+\frac{c_\sfc}{\beta_h}\right)\left(1+\frac{c_\sfa}{\alpha_h}\right)\left(1+\frac{c_\sfc}{\gamma_h}\right)\inf_{\tau_h\in W_h}\norm{\omega-\tau_h}_W.\nonumber
\end{gather}
\begin{proof}
The proof of this Theorem will be given in a series of Propositions \ref{prop:vorticityerrorbound} to \ref{prop:boundsigmas}.
\end{proof}
\end{theorem}
\begin{proposition}[\textbf{Vorticity error bound}]\label{prop:vorticityerrorbound}
Let $\sigma_h\in Z_{W_h}^\perp$, the error for vorticity is bounded by
\begin{equation}
\label{vorticityerrorbound}
\norm{\omega-\omega_h}_W\leq\left(1+\frac{c_\sfa}{\alpha_h}\right)\inf_{\sigma_h\in Z_{W_h}^\perp}\norm{\omega-\sigma_h}_W.
\end{equation}
\begin{proof}
Subtract the velocity-vorticity relation in the finite dimensional problem \eqref{discretedouble1} from that of the continuous problem \eqref{saddlepoint1}, we get
\[
\sfa(\tau_h,\omega-\omega_h)-\sfc(\tau_h,u-u_h)=0,\quad \forall\tau_h\in Z_{W_h}\subset Z_W.
\]
Bound $\sigma_h-\omega_h\in Z_{W_h}$ using inf-sup condition \eqref{infsupa}, we get
\[
\begin{aligned}
\alpha_h\norm{\sigma_h-\omega_h}_W&\leq\sup_{\tau_h\in Z_{W_h}}\frac{\sfa(\tau_h,\sigma_h-\omega_h)}{\norm{\tau_h}_W}\\
&=\sup_{\tau_h\in Z_{W_h}}\frac{\sfa(\tau_h,\sigma_h-\omega)+\sfa(\tau_h,\omega-\omega_h)}{\norm{\tau_h}_W}\\
&=\sup_{\tau_h\in Z_{W_h}}\frac{\sfa(\tau_h,\sigma_h-\omega)+\sfc(\tau_h,u-u_h)}{\norm{\tau_h}_W}.
\end{aligned}
\]
The last term vanishes since $\tau_h\in Z_{W_h}$ and $Z_{W_h}\subset Z_W$, hence $\alpha_h\norm{\sigma_h-\omega_h}_W\leq c_\sfa\norm{\omega-\sigma_h}_W$.
By the triangle inequality and the infimum over all $\sigma_h\in Z_{W_h}^\perp$ we obtain \eqref{vorticityerrorbound}.
\end{proof}
\end{proposition}
\begin{proposition}[\textbf{Velocity error bound}]\label{prop:velocityerrorbound}
Let $s_h\in Z_{V_h}^\perp$, the error for velocity is bounded by
\begin{equation}
\label{velocityerrorbound}
\norm{u-u_h}_V\leq\left(1+\frac{c_\sfc}{\gamma_h}\right)\inf_{s_h\in Z_{V_h}^\perp}\norm{u-s_h}_V+\frac{c_\sfa}{\gamma_h}\left(1+\frac{c_\sfa}{\alpha_h}\right)\inf_{\sigma_h\in Z_{W_h}^\perp}\norm{\omega-\sigma_h}_W.
\end{equation}
\begin{proof}
Use the inf-sup condition \eqref{infsupc} to bound $s_h-u_h\in Z_V$,
\[
\begin{aligned}
\gamma_h\norm{s_h-u_h}_V&\leq\sup_{\tau_h\in Z_{W_h}}\frac{\sfc(\tau_h,s_h-u_h)}{\norm{\tau_h}_W}\\
&=\sup_{\tau_h\in Z_{W_h}}\frac{\sfc(\tau_h,s_h-u)+\sfc(\tau_h,u-u_h)}{\norm{\tau_h}_W}\\
&=\sup_{\tau_h\in Z_{W_h}}\frac{\sfc(\tau_h,s_h-u)+\sfa(\tau_h,\omega-\omega_h)}{\norm{\tau_h}_W}\\
&\leq c_\sfc\norm{u-s_h}_V+c_\sfa\norm{\omega-\omega_h}_W.
\end{aligned}
\]
By triangle inequality, estimate \eqref{vorticityerrorbound} and the infimum over all $s_h\in Z_{V_h}^\perp$, we obtain \eqref{velocityerrorbound}.
\end{proof}
\end{proposition}
\begin{proposition}[\textbf{Pressure error bound}]\label{prop:pressureerrorbound}
The error for pressure is bounded by
\begin{equation}
\label{pressureerrorbound}
\begin{aligned}
\norm{p-p_h}_Q\leq&\left(1+\frac{c_\sfb}{\beta_h}\right)\inf_{q_h\in Q_h}\norm{p-q_h}_Q+\frac{c_\sfe}{\beta_h}\left(1+\frac{c_\sfc}{\gamma_h}\right)\inf_{s_h\in Z_{V_h}^\perp}\norm{u-s_h}_V\\
&+\left(\frac{c_\sfa}{\gamma_h}+\frac{c_\sfc}{\beta_h}\right)\left(1+\frac{c_\sfa}{\alpha_h}\right)\inf_{\sigma_h\in Z_{W_h}^\perp}\norm{\omega-\sigma_h}_W.
\end{aligned}
\end{equation}
\begin{proof}
Subtract \eqref{discretedouble2} from \eqref{saddlepoint2}, we get
\[
\sfc(\omega-\omega_h,v_h)+\sfe(v_h,u-u_h)+\sfb(v_h,p-p_h)=0,\quad\forall v_h\in V_h.
\]
So for $q_h\in Q_h$ we have
\[
\sfb(v_h,q_h-p_h)=-\sfc(\omega-\omega_h,v_h)-\sfe(v_h,u-u_h)-\sfb(v_h,p-q_h).
\]
Use this and the inf-sup condition \eqref{infsupb} to bound $q_h-p_h\in Q_h$,
\[
\begin{aligned}
\beta_h\norm{q_h-p_h}_Q&\leq\sup_{v_h\in V_h}\frac{\sfb(v_h,q_h-p_h)}{\norm{v_h}_V}\\
&=\sup_{v_h\in V_h}\frac{-\sfc(\omega-\omega_h,v_h)-\sfe(v_h,u-u_h)-\sfb(v_h,p-q_h)}{\norm{v_h}_V}\\
&\leq c_\sfc\norm{\omega-\omega_h}_W+c_\sfe\norm{u-u_h}_V+c_\sfb\norm{p-q_h}_Q.
\end{aligned}
\]
By triangle inequality, estimates \eqref{vorticityerrorbound} and \eqref{velocityerrorbound}, and the infimum over all $q_h\in Q_h$, we obtain \eqref{pressureerrorbound}.
\end{proof}
\end{proposition}
Next we replace the infimums over $\sigma_h\in Z_{W_h}^\perp$ and $s_h\in Z_{V_h}^\perp$ by best approximation errors.
\begin{proposition}\label{prop:boundsigmas}
The terms $\inf_{\sigma_h\in Z_{W_h}^\perp}\norm{\omega-\sigma_h}_W$ and $\inf_{s_h\in Z_{V_h}^\perp}\norm{u-s_h}_V$ are bounded by the best approximation estimates $\inf_{\tau_h\in W_h}\norm{\omega-\tau_h}_W$ and $\inf_{v_h\in V_h}\norm{u-v_h}_V$, using the inf-sup conditions \eqref{infsupb} and \eqref{infsupc}, as
\begin{equation}
\label{boundsigma}
\inf_{\sigma_h\in Z_{W_h}^\perp}\norm{\omega-\sigma_h}_W\leq\left(1+\frac{c_\sfc}{\gamma_h}\right)\inf_{\tau_h\in W_h}\norm{\omega-\tau_h}_W,
\end{equation}
\begin{equation}
\label{bounds}
\inf_{s_h\in Z_{V_h}^\perp}\norm{u-s_h}_V\leq\left(1+\frac{c_\sfb}{\beta_h}\right)\inf_{v_h\in V_h}\norm{u-v_h}_V.
\end{equation}
\begin{proof}
Take $\tau_h\in W_h$, then there exists a $\kappa_h\in W_h$ such that
\[
\sfc(\kappa_h,v_{\mathcal{Z}_h})=\sfc(\omega-\tau_h,v_{\mathcal{Z}_h}),\quad\forall v_{\mathcal{Z}_h}\in Z_{V_h}.
\]
This is equivalent to
\[
\sfc(\kappa_h+\tau_h,v_{\mathcal{Z}_h})=\sfc(\omega,v_{\mathcal{Z}_h})=\sff(v_{\mathcal{Z}_h}),\quad\forall v_{\mathcal{Z}_h}\in Z_{V_h},
\]
which shows that $\sigma_h=\kappa_h+\tau_h\in Z_{W_h}^\perp$. 
We can bound $\norm{\kappa_h}_W$ using the discrete inf-sup condition as follows
\[
\gamma_h\norm{\kappa_h}_W\leq\sup_{v\in Z_{V_h}}\frac{\sfc(\kappa_h,v)}{\norm{v}_V}=\sup_{v\in Z_{V_h}}\frac{\sfc(\omega-\tau_h,v)}{\norm{v}_V}\leq c_\sfc\norm{\omega-\tau_h}_V.
\]
By triangle inequality and since $\tau_h\in W_h$ was arbitrary, we find \eqref{boundsigma}. A similar proof holds for \eqref{bounds} (see also \cite{brezzifortin}, Proposition 2.5).
\end{proof}
\end{proposition}

Additionally, following section 7.7.6 in \cite{bochevgunzburger}, we have the following $L^2\Lambda^{k}(\Omega)$ error estimates for the curl of vorticity and divergence of velocity,
\begin{proposition}
The errors of the curl of vorticity and divergence of velocity are bounded by their best approximation estimates,
\begin{align}
\norm{\ud(\omega-\omega_h)}_{L^2\Lambda^{n-1}}&\leq\inf_{\tau_h\in W_h}\norm{\ud(\omega-\tau_h)}_{L^2\Lambda^{n-1}},\label{curlvorticity}\\
\norm{\ud(u-u_h)}_{L^2\Lambda^{n}}&\leq\inf_{v_h\in V_h}\norm{\ud(u-v_h)}_{L^2\Lambda^{n}}.\label{divvelocity}
\end{align}
\begin{proof}
Choose $v=v_{\mathcal{Z}_h}\in Z_{V_h}$ in \eqref{saddlepoint2} and \eqref{discretedouble2} and subtract these. Set $v_{\mathcal{Z}_h}=\ud\tau_h$, this gives the orthogonality relation
\[
(\ud(\omega-\omega_h),\ud\tau_h))_\Omega=0,\quad\forall\tau_h\in W_h.
\]
Substitute this into the following Cauchy-Schwarz inequality,
\[
\begin{aligned}
\norm{\ud(\omega-\omega_h)}^2_{L^2\Lambda^{n-1}}&=(\ud(\omega-\omega_h),\ud(\omega-\tau_h))_\Omega\\
&\leq\norm{\ud(\omega-\omega_h)}_{L^2\Lambda^{n-1}}\norm{\ud(\omega-\tau_h)}_{L^2\Lambda^{n-1}},\quad\forall\tau_h\in W_h,
\end{aligned}
\]
and \eqref{curlvorticity} follows. Next choose $q_h=\ud v_h\in Q_h$ in \eqref{saddlepoint3} and \eqref{discretedouble3} and subtract these. Then \eqref{divvelocity} follows again from the Cauchy-Schwarz inequality.
\end{proof}
\end{proposition}


Because the projections of respectively $\omega,\ u,$ and $p$, belong to the finite dimensional subspaces $W_h\subset W$, $V_h\subset V$ and $Q_h\subset Q$, the best approximation errors can be bounded using the interpolation errors,
\[
\inf_{\tau_h\in W_h}\norm{\omega-\tau_h}_W\leq\norm{\omega-\pi_h\omega}_W,\ \inf_{v_h\in V_h}\norm{u-v_h}_V\leq\norm{u-\pi_h u}_V,\ \inf_{q_h\in Q_h}\norm{p-q_h}_Q\leq\norm{p-\pi_hp}_Q,
\]
and therefore we obtain the following optimal a priori error estimates,
\begin{equation}
\norm{\omega-\omega_h}_W=\mathcal{O}(h^{N}),\quad \norm{u-u_h}_V=\mathcal{O}(h^N), \quad \norm{p-p_h}_Q=\mathcal{O}(h^N).
\end{equation}
So the convergence rates for the approximations are equal to those of the interpolations, \eqref{L2interpolation}, \eqref{Hinterpolation}, thus we obtained optimal convergence. The error estimates were obtained independent of the chosen types of boundary conditions.

\begin{remark}
In contrast to \cite{arnold2011,dubois2003b} and \cite{bochevgunzburger}, Table 7.5, where Raviart-Thomas elements were used, the proposed compatible method has provably optimal convergence also with standard velocity boundary conditions and with non-affine mappings.
\end{remark}

\section{Numerical Results}\label{sec:numericalresults}
Now that the compatible spectral discretization method and its a priori error estimates are derived, we perform a series of test problems to show optimal convergence behavior. Purpose of the testcases is to show convergence behavior in case of various boundary conditions and in case of curvilinear meshes. In all cases we show optimal convergence.

The first three testcases originate from a recent paper by Arnold et al \cite{arnold2011}, where suboptimal convergence is shown for normal velocity - tangential boundary conditions in vector Poisson and Stokes problems, when using Raviart-Thomas elements \cite{raviartthomas1977}. Since Raviart-Thomas elements are the most popular $H(\mathrm{div,\Omega})$ conforming elements, we compare our method to these results.

\subsection{Vector Poisson problems}
\begin{figure}[htbp]
\centering
\includegraphics[width=.8\textwidth]{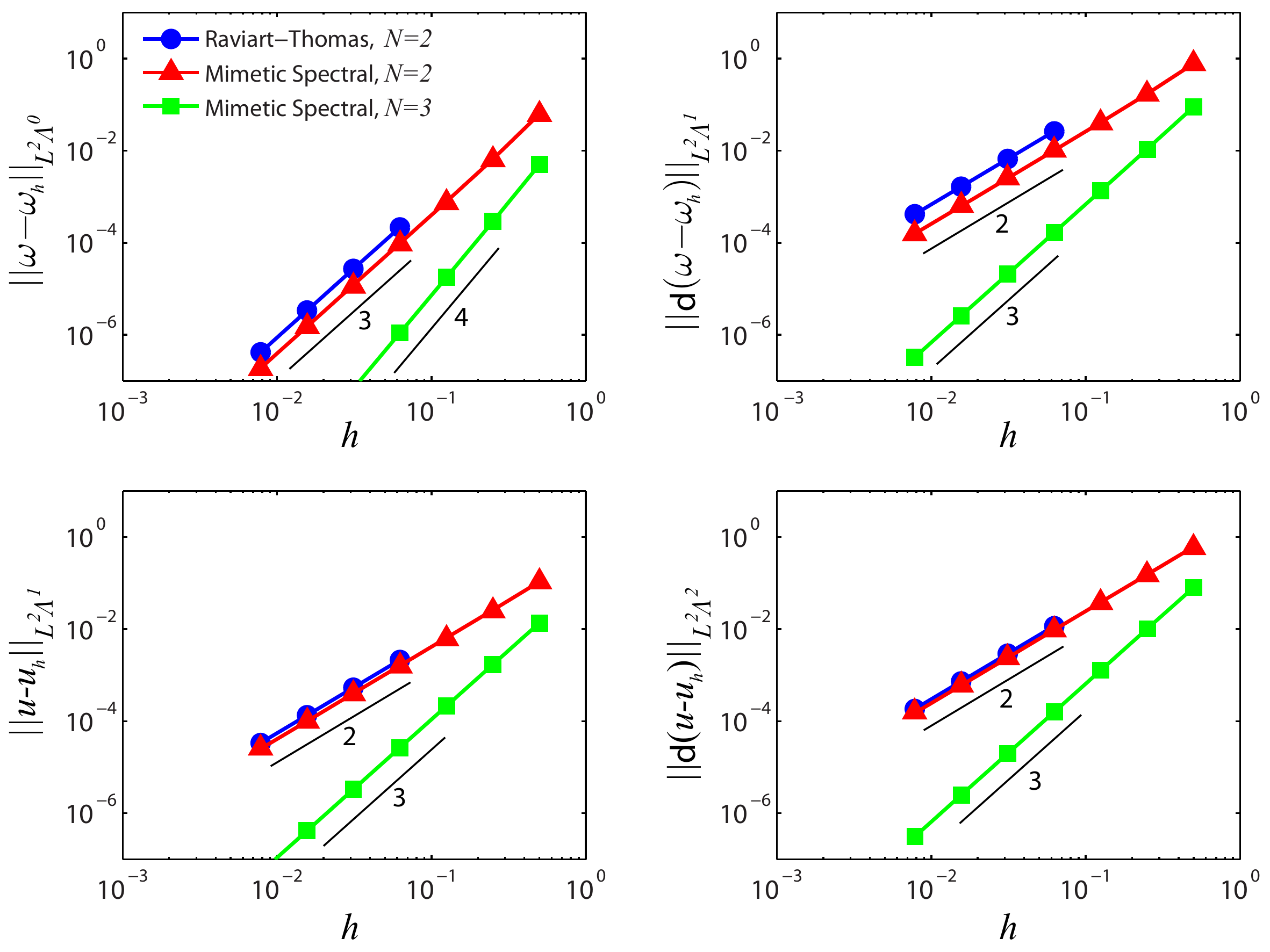}
\caption{Comparison of the $h$-convergence between Raviart-Thomas and Mimetic spectral elements for the 2D 1-form Poisson problem with tangential velocity - divergence-free boundary conditions.}
\label{fig:afgcase1}
\end{figure}

\figref{fig:afgcase1} shows the result of the vector Poisson problem \eqref{poisson} on $\Omega=[0,1]^2$ with coordinates $\mathbf{x}:=(x,y)$, for a 1-form $u\in H\Lambda^1(\Omega)$, where $\Gamma=\Gamma_2$, i.e. with tangential velocity - divergence-free boundary conditions ($\tr\star u=0,\ \tr\star\ud u=0$). The corresponding solution is given by
\begin{align}
u^{(1)}&=-v(\mathbf{x})\,\ud x+u(\mathbf{x})\,\ud y\nonumber\\
&=-(2\sin\pi x\cos\pi y)\,\ud x+ (\cos\pi x\sin\pi y)\,\ud y.
\end{align}
Both Raviart-Thomas and mimetic spectral element methods show optimal convergence rates. All results of this and the following two problems where obtained on the same quadrilateral mesh of $2^n\times 2^n$ subsquares, $n=1,2,3,4,\hdots$ just like the reference solutions from \cite{arnold2011}.

\begin{figure}[htbp]
\centering
\includegraphics[width=.8\textwidth]{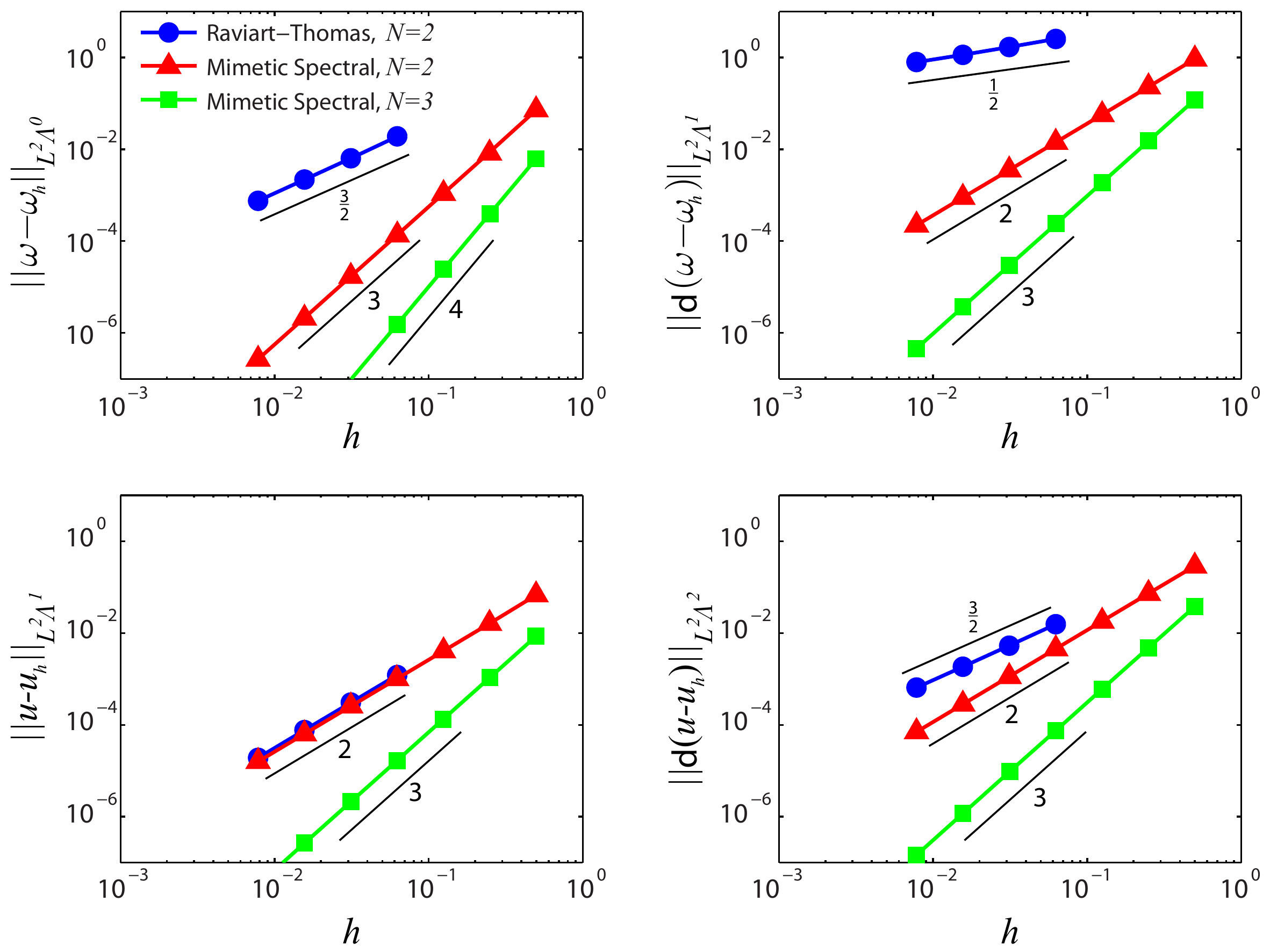}
\caption{Comparison of the $h$-convergence between Raviart-Thomas and Mimetic spectral elements for the 2D 1-form Poisson problem with tangential velocity - normal velocity boundary conditions.}
\label{fig:afgcase2}
\end{figure}

\figref{fig:afgcase2} shows again results for the vector Poisson problem for a 1-form, but now in combination with normal velocity - tangential velocity boundary conditions ($\tr u=0,\ \tr\star u=0$), so $\Gamma=\Gamma_1$. The corresponding manufactured solution is
\begin{align}
u^{(1)}&=-v(\mathbf{x})\,\ud x+u(\mathbf{x})\,\ud y\nonumber\\
&=-(\sin\pi x\sin\pi y)\,\ud x+ (\sin\pi x\sin\pi y)\,\ud y.
\end{align}
The compatible spectral discretization method again shows optimal convergence, as was expected from the above analysis. The Raviart-Thomas elements only show suboptimal convergence in case of velocity boundary conditions. This suboptimality was proven in \cite{arnold2011}. Especially for $\omega_h$ and $\ud\omega_h$ the current method outperforms the Raviart-Thomas elements, with a difference in rate of convergence of $\tfrac{3}{2}$.

\subsection{Stokes problems}
\begin{figure}[htbp]
\centering
\includegraphics[width=.8\textwidth]{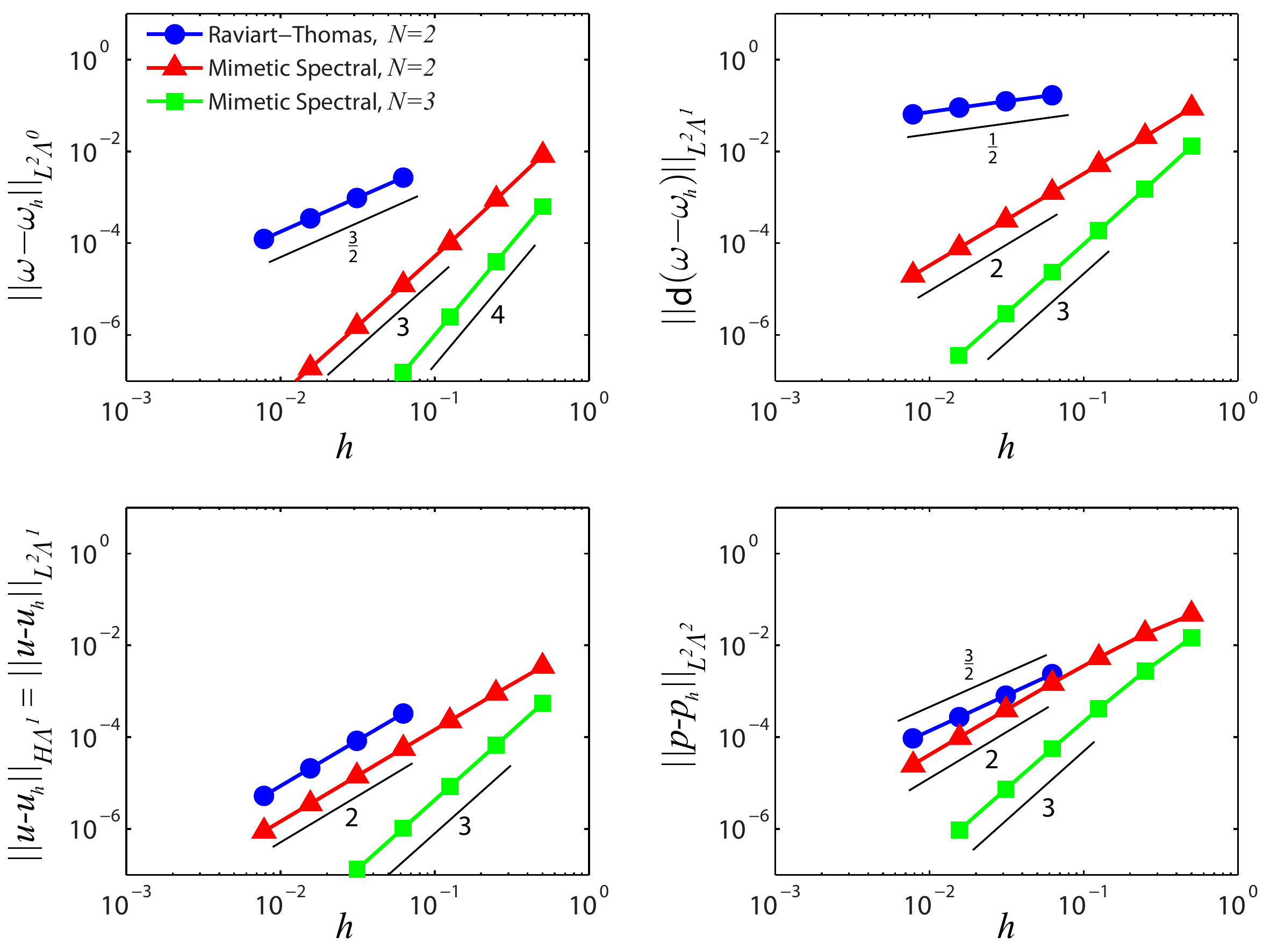}
\caption{Comparison of the $h$-convergence between Raviart-Thomas and Mimetic spectral element projections for the 2D Stokes problem with normal velocity - tangential velocity boundary conditions.}
\label{fig:afgcase3}
\end{figure}

The same difference in convergence behavior is found for the Stokes problem, where $\Gamma=\Gamma_1$, i.e. with normal velocity - tangential velocity boundary conditions, see \figref{fig:afgcase3}. Again $\Omega$ is the unit square, and the velocity and pressure fields are given by
\begin{align}
u^{(1)}&=-v(\mathbf{x})\,\ud x+u(\mathbf{x})\,\ud y\nonumber\\
&=-\left(2y^2(y-1)^2x(2x-1)(x-1)\right)\,\ud x+ \left(-2x^2(x-1)^2y(2y-1)(y-1)\right)\,\ud y,\\
p^{(2)}&=p(\mathbf{x})\,\ud x\wedge\ud y=\left((x-\tfrac{1}{2})^5+(y-\tfrac{1}{2})^5\right)\,\ud x\wedge\ud y.
\end{align}
While for velocity both methods show optimal convergence, for pressure a difference of $\tfrac{1}{2}$ is noticed in the rate of convergence and for vorticity and the curl of vorticity again a difference in rate of convergence of $\tfrac{3}{2}$ is revealed.
The error in divergence of velocity is not shown here for the Stokes problem, because the method is pointwise divergence-free up to machine precision. Special attention to this property is given in \cite{kreeft2012}.

We would like to remark is that the results shown in \figref{fig:afgcase3} are independent of the kind of boundary conditions used. Table~\ref{tab:bcresults} shows the results of vorticity for all types of admissible boundary conditions.

\begin{table}[htbp]
\centering
\begin{tabular}{c|c|c|c||c}
 normal velocity & tangential velocity & vorticity & \ \ vorticity\ \  & convergence\\
 tangential velocity & pressure & normal velocity & pressure & rate\\ \hline\hline
   1.0280e-04 & 1.0109e-04 & 1.0030e-04 & 1.0035e-04 & 3.14 \\
   1.2445e-05 & 1.2410e-05 & 1.2364e-05 & 1.2375e-05 & 3.05 \\
   1.5424e-06 & 1.5426e-06 & 1.5399e-06 & 1.5416e-06 & 3.01 \\
   1.9238e-07 & 1.9247e-07 & 1.9230e-07 & 1.9255e-07 & 3.00 \\
   2.4035e-08 & 2.4042e-08 & 2.4032e-08 & 2.4065e-08 & 3.00 \\
\end{tabular}
\caption{This table shows the vorticity error $\norm{\omega-\omega_h}_{L^2\Lambda^{0}}$ obtained using the four types of boundary conditions given in \eqref{boundaryconditions}. The results are obtained on an uniform Cartesian mesh with $N=2$ and $h=\tfrac{1}{8},\tfrac{1}{16},\tfrac{1}{32},\tfrac{1}{64},\tfrac{1}{128}$. All four cases show third order convergence.}
\label{tab:bcresults}
\end{table}

\begin{figure}[htbp]
\centering
\includegraphics[width=0.8\textwidth]{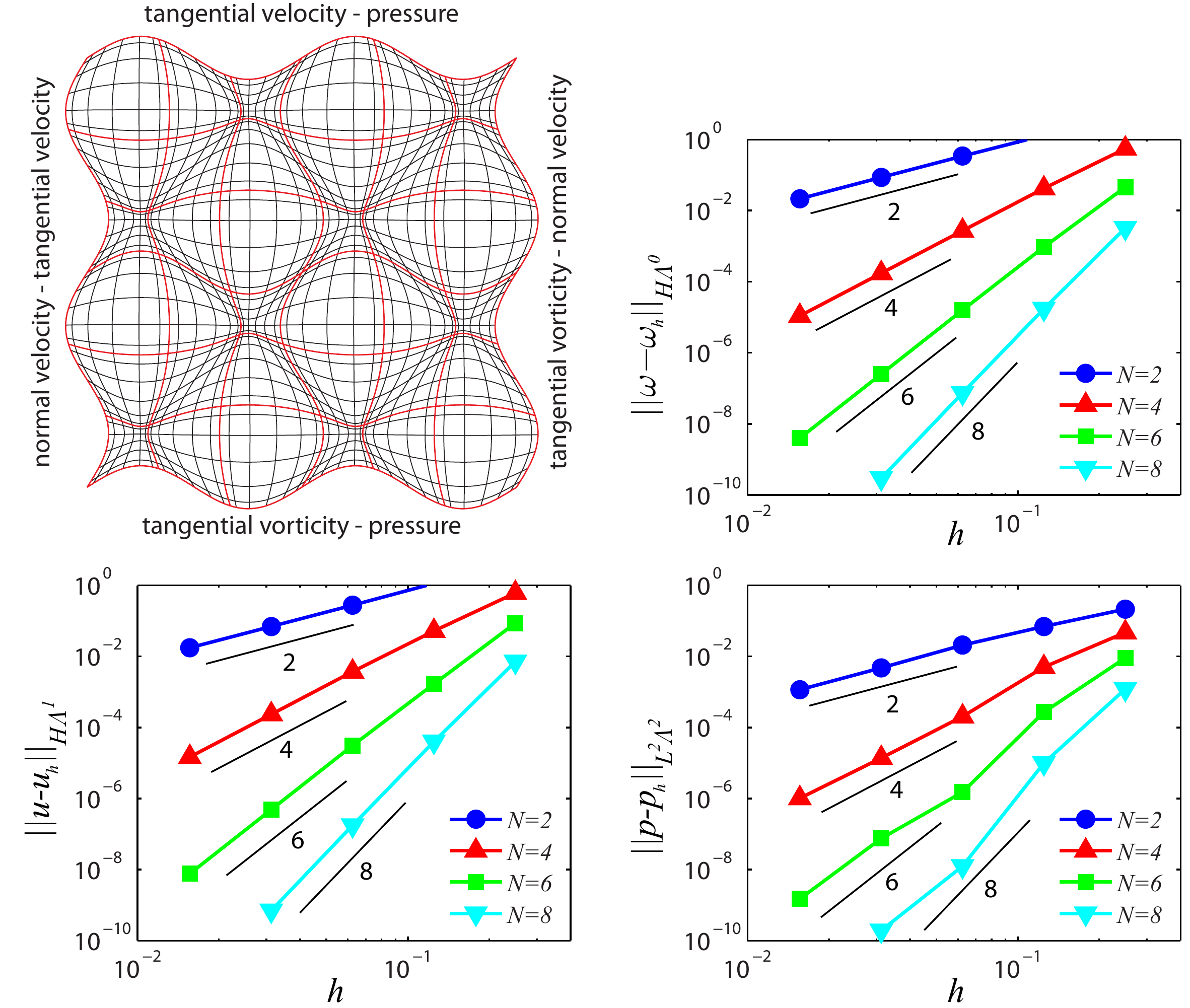}
\caption{Upper left figure show the computational domain with a $6\times 6$ element mesh of $N=6$. Furthermore the velocity, vorticity and pressure $h$-convergence results are shown of Stokes problem \eqref{testcase4}. All variables are tested on grids with $N=2,4,6$ and 8.}
\label{fig:allinonestokes}
\end{figure}

The next testcase reveals the optimal convergence in case of higher-order approximation on curvilinear quadrilateral meshes for all admissible types of boundary conditions. The manufactured solution Stokes problem is given on a curvilinear domain, defined by the mapping $(x,y)=\Phi(\xi,\eta)$,
\begin{subequations}
\begin{align}
x(\xi,\eta)&=\tfrac{1}{2}+\tfrac{1}{2}\left(\xi+\tfrac{1}{10}\cos(2\pi\xi)\sin(2\pi\eta)\right),\\
y(\xi,\eta)&=\tfrac{1}{2}+\tfrac{1}{2}\left(\eta+\tfrac{1}{10}\sin(2\pi\xi)\cos(2\pi\eta)\right).
\end{align}
\end{subequations}
A $6\times6$ element $N=6$ mesh is show in \figref{fig:allinonestokes}. Each side of the domain has a different type of boundary condition, so $\Gamma=\Gamma_1\cup\Gamma_2\cup\Gamma_3\cup\Gamma_4$, as shown in the same figure and listed in \eqref{boundaryconditions}.
The solutions of vorticity $\omega\in\Lambda^0(\Omega)$, velocity $u\in\Lambda^1(\Omega)$ and pressure $p\in\Lambda^2(\Omega)$ are given by
\begin{subequations}
\label{testcase4}
\begin{align}
\omega^{(0)}&=\tfrac{3}{2}\pi\sin(\tfrac{3}{2}\pi x)\sin(\tfrac{3}{2}\pi y),\\
u^{(1)}&=-\left(\cos(\tfrac{3}{2}\pi x)\sin(\tfrac{3}{2}\pi y)\right)\,\ud x+\left(2\sin(\tfrac{3}{2}\pi x)\cos(\tfrac{3}{2}\pi y)\right)\,\ud y,\\
p^{(2)}&=\left(\sin(\pi x)\sin(\pi y)\right)\,\ud x\wedge\ud y.
\end{align}
\end{subequations}
They lead to nonzero body force $f\in\Lambda^1(\Omega)$ and mass source $g\in\Lambda^2(\Omega)$. \figref{fig:allinonestokes} shows the convergence of the vorticity $\omega_h\in\Lambda_h^{0}(\Omega;C_0)$, velocity $u_h\in\Lambda^{1}_h(\Omega;C_1)$ and pressure $p_h\in\Lambda^2_h(\Omega;C_2)$. The errors for the vorticity and velocity are measured in the $H\Lambda^k$-norm, i.e. $\Vert\omega-\omega_h\Vert_{H\Lambda^{0}}$, and $\Vert u-u_h\Vert_{H\Lambda^{1}}$, respectively, and the error of the pressure is given in the $L^2\Lambda^2$-norm. In \figref{fig:allinonestokes} convergence rates are added which show the \emph{optimal} $h$-convergence behavior of the Stokes problem on a curvilinear domain with curvilinear grid and all four types of boundary conditions.

\section{Concluding remark}
Optimal approximation of the Stokes problem for all admissible boundary conditions essentially hinges on the construction of a conforming discrete Hodge decomposition, $\Lambda^k_h=\mathcal{Z}^k_h\oplus\mathcal{Z}^{k,\perp}_h$ and a discrete Poicar\'{e} inequality, that are based on the bijection of the exterior derivative on the conforming subspace, $\ud:\mathcal{Z}_h^{k,\perp}\rightarrow\mathcal{B}_h^{k+1}$. Ensuring these properties result in a compatible discretization method, and relied on the construction of a bounded projection operator, $\pi_h:\Lambda^k(\Omega)\rightarrow\Lambda_h^k(\Omega;C_k)$, that commutes with the exterior derivative, $\pi_h\ud=\reconstruction\delta\reduction=\ud\pi_h$. So the compatibility is based on the bijection of the coboundary operator, $\delta:Z^{k,\perp}\rightarrow B^{k+1}$, and the construction of interpolatory basis functions.
From this it follows that, $\mathcal{B}_h^{k+1}\subset\mathcal{B}^{k+1}$, $\mathcal{Z}_h^k\subset\mathcal{Z}^k$ and $\mathcal{Z}_h^{k,\perp}\subset\mathcal{Z}^{k,\perp}$. From these properties the rest follows.

For piecewise sufficiently smooth mappings, the optimal conference rates hold on curvilinear grids as well, since the pullback operator of the map from a curvilinear domain to the Cartesian frame commutes with the projection operator. Any projection (discretization) with these properties will yield similar results as described in this paper.

\section*{Acknowledgments}
We would like to thank Pavel Bochev for the fruitful discussions on mimetic schemes, boundary conditions and error estimates.

\bibliographystyle{abbrv}
\bibliography{./literature}

\end{document}